\documentclass[a4paper,twoside]{article}
\usepackage{a4}
\usepackage{amssymb}
\usepackage{amsmath}
\usepackage{upref}
\usepackage[active]{srcltx}
\allowdisplaybreaks[2] 
\usepackage[colorlinks,citecolor=blue,linkcolor=blue]{hyperref}
\newcount\minutes \newcount\hours
\hours=\time
\divide\hours 60
\minutes=\hours
\multiply\minutes -60
\advance\minutes \time
\newcommand{\klockan}{\the\hours:{\ifnum\minutes<10 0\fi}\the\minutes}
\newcommand{\tid}{\today\ \klockan}
\newcommand{\prtid}{\smash{\raise 10mm \hbox{\LaTeX ed \tid}}}
\renewcommand{\prtid}{}
\makeatletter
\def\sectionmark#1{} 
\def\subsectionmark#1{}
\newcommand{\sectnr}{\ifnum \c@secnumdepth >\z@
                 \thesection.\hskip 1em\relax \fi}
\def\@evenhead{\footnotesize\rm\thepage\hfil\leftmark\hfil\llap{\prtid}}
\def\@oddhead{\footnotesize\rm\rlap{\prtid}\hfil\rightmark\hfil\thepage}
\def\tableofcontents{\section*{Contents} 
 \@starttoc{toc}}
\makeatother
\makeatletter
\def\@biblabel#1{#1.}
\makeatother
\makeatletter
\let\Thebibliography=\thebibliography
\renewcommand{\thebibliography}[1]{\def\@mkboth##1##2{}\Thebibliography{#1}
\addcontentsline{toc}{section}{References}
\frenchspacing 
\setlength{\@topsep}{0pt}
\setlength{\itemsep}{0pt}%
\setlength{\parskip}{0pt plus 2pt}%
}
\makeatother
\makeatletter
\def\mdots@{\mathinner.\nonscript\!.%
 \ifx\next,.\else\ifx\next;.\else\ifx\next..\else
 \nonscript\!\mathinner.\fi\fi\fi}
\let\ldots\mdots@
\let\cdots\mdots@
\let\dotso\mdots@
\let\dotsb\mdots@
\let\dotsm\mdots@
\let\dotsc\mdots@
\def\vdots{\vbox{\baselineskip2.8\p@ \lineskiplimit\z@
    \kern6\p@\hbox{.}\hbox{.}\hbox{.}\kern3\p@}}
\def\ddots{\mathinner{\mkern1mu\raise8.6\p@\vbox{\kern7\p@\hbox{.}}%
    \raise5.8\p@\hbox{.}\raise3\p@\hbox{.}\mkern1mu}}
\makeatother
\makeatletter
\let\Enumerate=\enumerate
\renewcommand{\enumerate}{\Enumerate%
\setlength{\@topsep}{0pt}
\setlength{\itemsep}{0pt}%
\setlength{\parskip}{0pt plus 1pt}%
\renewcommand{\theenumi}{\textup{(\alph{enumi})}}%
\renewcommand{\labelenumi}{\theenumi}%
}
\let\endEnumerate=\endenumerate
\renewcommand{\endenumerate}{\endEnumerate\unskip}
\makeatother
\makeatletter
\def\@seccntformat#1{\csname the#1\endcsname.\quad}
\makeatother
\newcommand{\authortitle}[2]{\author{#1}\title{#2}\markboth{#1}{#2}}
\newcommand{\art}[6]{{\sc #1, \rm #2, \it #3 \bf #4 \rm (#5), \mbox{#6}.}}
\newcommand{\artin}[3]{{\sc #1, \rm #2,  in #3.}}
\newcommand{\artnopt}[6]{{\sc #1, \rm #2, \it #3 \bf #4 \rm (#5), \mbox{#6}}}
\newcommand{\arttoappear}[3]{{\sc #1, \rm #2, to appear in \it #3}}
\newcommand{\auth}[2]{{#1, #2.}}
\newcommand{\book}[3]{{\sc #1, \it #2, \rm #3.}}
\newcommand{\AND}{{\rm and }}
\RequirePackage{amsthm}
\newtheoremstyle{descriptive}%
  {\topsep}   
  {\topsep}   
  {\rmfamily} 
  {}          
  {\bfseries} 
  {.}         
  { }         
  {}          
\newtheoremstyle{propositional}%
  {\topsep}   
  {\topsep}   
  {\itshape}  
  {}          
  {\bfseries} 
  {.}         
  { }         
  {}          
\theoremstyle{propositional}
\newtheorem{thm}{Theorem}[section]
\newtheorem{prop}[thm]{Proposition}
\newtheorem{lem}[thm]{Lemma}
\newtheorem{cor}[thm]{Corollary}
\theoremstyle{descriptive}
\newtheorem{deff}[thm]{Definition}
\newtheorem{example}[thm]{Example}
\makeatletter
\renewenvironment{proof}[1][\proofname]{\par
  \pushQED{\qed}%
  \normalfont 
  \trivlist
  \item[\hskip\labelsep
        \itshape
    #1\@addpunct{.}]\ignorespaces
}{%
  \popQED\endtrivlist\@endpefalse
}
\makeatother
\newcommand{\setm}{\setminus}
\renewcommand{\subsetneq}{\varsubsetneq}
\renewcommand{\emptyset}{\varnothing}
\def\vint{\mathop{\mathchoice%
          {\setbox0\hbox{$\displaystyle\intop$}\kern 0.22\wd0%
           \vcenter{\hrule width 0.6\wd0}\kern -0.82\wd0}%
          {\setbox0\hbox{$\textstyle\intop$}\kern 0.2\wd0%
           \vcenter{\hrule width 0.6\wd0}\kern -0.8\wd0}%
          {\setbox0\hbox{$\scriptstyle\intop$}\kern 0.2\wd0%
           \vcenter{\hrule width 0.6\wd0}\kern -0.8\wd0}%
          {\setbox0\hbox{$\scriptscriptstyle\intop$}\kern 0.2\wd0%
           \vcenter{\hrule width 0.6\wd0}\kern -0.8\wd0}}%
          \mathopen{}\int}
\newcommand{\Cp}{{C_p}}
\newcommand{\bCp}{{\protect\itoverline{C}_p}}
\DeclareMathOperator{\diam}{diam}
\DeclareMathOperator{\Div}{div}
\DeclareMathOperator{\Irr}{Irr}
\DeclareMathOperator{\res}{res}
\DeclareMathOperator*{\essliminf}{ess\,lim\,inf}
\DeclareMathOperator*{\essinf}{ess\,inf}
\newcommand{\bdy}{\partial}
\newcommand{\loc}{_{\rm loc}}
\DeclareMathOperator{\id}{id}
{\catcode`p =12 \catcode`t =12 \gdef\eeaa#1pt{#1}}      
\def\accentadjtext#1{\setbox0\hbox{$#1$}\kern   
                \expandafter\eeaa\the\fontdimen1\textfont1 \ht0 }
\def\accentadjscript#1{\setbox0\hbox{$#1$}\kern 
                \expandafter\eeaa\the\fontdimen1\scriptfont1 \ht0 }
\def\accentadjscriptscript#1{\setbox0\hbox{$#1$}\kern   
                \expandafter\eeaa\the\fontdimen1\scriptscriptfont1 \ht0 }
\def\accentadjtextback#1{\setbox0\hbox{$#1$}\kern       
                -\expandafter\eeaa\the\fontdimen1\textfont1 \ht0 }
\def\accentadjscriptback#1{\setbox0\hbox{$#1$}\kern     
                -\expandafter\eeaa\the\fontdimen1\scriptfont1 \ht0 }
\def\accentadjscriptscriptback#1{\setbox0\hbox{$#1$}\kern 
                -\expandafter\eeaa\the\fontdimen1\scriptscriptfont1 \ht0 }
\def\itoverline#1{{\mathsurround0pt\mathchoice
        {\rlap{$\accentadjtext{\displaystyle #1}
                \accentadjtext{\vrule height1.593pt}
                \overline{\phantom{\displaystyle #1}
                \accentadjtextback{\displaystyle #1}}$}{#1}}
        {\rlap{$\accentadjtext{\textstyle #1}
                \accentadjtext{\vrule height1.593pt}
                \overline{\phantom{\textstyle #1}
                \accentadjtextback{\textstyle #1}}$}{#1}}
        {\rlap{$\accentadjscript{\scriptstyle #1}
                \accentadjscript{\vrule height1.593pt}
                \overline{\phantom{\scriptstyle #1}
                \accentadjscriptback{\scriptstyle #1}}$}{#1}}
        {\rlap{$\accentadjscriptscript{\scriptscriptstyle #1}
                \accentadjscriptscript{\vrule height1.593pt}
                \overline{\phantom{\scriptscriptstyle #1}
                \accentadjscriptscriptback{\scriptscriptstyle #1}}$}{#1}}}}
\def\itunderline#1{{\mathsurround0pt\mathchoice
        {\rlap{$\underline{\phantom{\displaystyle #1}
                \accentadjtextback{\displaystyle #1}}$}{#1}}
        {\rlap{$\underline{\phantom{\textstyle #1}
                \accentadjtextback{\textstyle #1}}$}{#1}}
        {\rlap{$\underline{\phantom{\scriptstyle #1}
                \accentadjscriptback{\scriptstyle #1}}$}{#1}}
        {\rlap{$\underline{\phantom{\scriptscriptstyle #1}
                \accentadjscriptscriptback{\scriptscriptstyle #1}}$}{#1}}}}
\newcommand{\limplus}{{\mathchoice{\vcenter{\hbox{$\scriptstyle +$}}}
  {\vcenter{\hbox{$\scriptstyle +$}}}
  {\vcenter{\hbox{$\scriptscriptstyle +$}}}
  {\vcenter{\hbox{$\scriptscriptstyle +$}}}
}}
\def\cprime{{\mathsurround0pt$'$}}
\newcommand{\ga}{\gamma}
\newcommand{\de}{\delta}
\newcommand{\eps}{\varepsilon}
\newcommand{\la}{\lambda}
\newcommand{\Om}{\Omega}
\newcommand{\clOmprim}{{\overline{\Om}\mspace{1mu}}'}
\newcommand{\clOmm}{{\overline{\Om}\mspace{1mu}}^M}
\newcommand{\clOmj}{{\overline{\Om}^j}}
\newcommand{\clOm}{{\overline{\Om}}}
\newcommand{\Omm}{{\Om^M}}
\newcommand{\bdym}{{\bdy^M}}
\newcommand{\bdyj}{{\bdy^j}}
\renewcommand{\phi}{\varphi}
\newcommand{\p}{{$p\mspace{1mu}$}}   
\newcommand{\R}{\mathbf{R}}
\newcommand{\eR}{{\overline{\R}}}
\newcommand{\Ga}{\Gamma}
\newcommand{\Np}{N^{1,p}}
\newcommand{\Dp}{D^p}
\newcommand{\Nploc}{N^{1,p}\loc}
\newcommand{\Lploc}{L^p\loc}
\newcommand{\Hpind}[1]{P_{#1}}      
\newcommand{\uHpind}[1]{\itoverline{P}_{#1}}      
\newcommand{\lHpind}[1]{\itunderline{P}_{#1}}      
\newcommand{\uP}{\itoverline{P}}     
\newcommand{\lP}{\itunderline{P}} 
\newcommand{\A}{\mathcal{A}}
\newcommand{\UU}{\mathcal{U}}
\newcommand{\ft}{\tilde{f}}
\newcommand{\dMto}{\overset{d_M}\longrightarrow}
\newcommand{\dM}{d_M}

\newcommand{\bdyone}{{\bdy^1}}
\newcommand{\tauone}{{\tau^1}}
\newcommand{\tauj}{{\tau^j}}
\newcommand{\toone}{{\overset{\tau^1}\longrightarrow}}
\newcommand{\clOmone}{{\overline{\Om}^1}}
\newcommand{\clGone}{{\itoverline{G}^1}}
\newcommand{\Omone}{{\Om^1}}
\newcommand{\Uone}{{U^1}}
\newcommand{\Utwo}{{U^2}}
\newcommand{\Vtwo}{{V^2}}
\newcommand{\clUjtwo}{{\overline{U}^2_j}}
\newcommand{\clOmQ}{{\overline{\Om}\mspace{1mu}}^Q}

\newcommand{\bdytwo}{{\bdy^2}}
\newcommand{\tautwo}{{\tau^2}}
\newcommand{\totwo}{{\overset{\tau^2}\longrightarrow}}
\newcommand{\clOmtwo}{{\overline{\Om}^2}}
\newcommand{\clOmptwo}{{\overline{\Om'}^2}}
\newcommand{\clVptwo}{{\overline{V'}^2}}
\newcommand{\clUtwo}{{\overline{U}^2}}

\newcommand{\Omtwo}{{\Om^2}}
\newcommand{\Omptwo}{{(\Om')^2}}
\newcommand{\Irrres}{\Irr_{\res}}

\newcommand{\uSpind}[1]{\itoverline{S}_{#1}}      
\newcommand{\lSpind}[1]{\itunderline{S}_{#1}}      

\newcommand{\xh}{\hat{x}}
\newcommand{\fh}{\hat{f}}
\newcommand{\uh}{\hat{u}}
\newcommand{\Gh}{\widehat{G}}
\newcommand{\Vh}{\widehat{V}}
\newcommand{\Kh}{\widehat{K}}
\newcommand{\DU}{\mathcal{DU}}%
\numberwithin{equation}{section}
\newcommand{\eqv}{\ensuremath{
\mathchoice{\quad \Longleftrightarrow \quad}{\Leftrightarrow}
                {\Leftrightarrow}{\Leftrightarrow}} }
\newcommand{\imp}{\ensuremath{\Rightarrow} }

\newenvironment{ack}{\medskip{\it Acknowledgement.}}{}

\begin{document}

\authortitle{Anders Bj\"orn}
{The Kellogg property  for \p-harmonic functions
with respect to the Mazurkiewicz boundary}
\title{The Kellogg property and boundary regularity for \p-harmonic functions
with respect to the Mazurkiewicz boundary and other compactifications}
\author{
Anders Bj\"orn \\
\it\small Department of Mathematics, Link\"oping University, \\
\it\small SE-581 83 Link\"oping, Sweden\/{\rm ;}
\it \small anders.bjorn@liu.se
}

\date{}
\maketitle

\noindent{\small
{\bf Abstract.} 
In this paper boundary regularity for \p-harmonic functions
is studied with respect to the Mazurkiewicz boundary and other compactifications.
In particular, the Kellogg property 
(which says that the set of irregular boundary points has capacity zero)
is obtained for a large class of compactifications,
but also two examples when it fails are given.
This study is done for complete metric spaces equipped with doubling measures
supporting a \p-Poincar\'e inequality,
but the results are new also in unweighted Euclidean spaces.
} 

\bigskip
\noindent
    {\small \emph{Key words and phrases}:
      boundary regularity,
      capacity,
Dirichlet problem,
doubling measure,
Kellogg property,
metric space,
\p-harmonic function,
Poincar\'e inequality,
regular point,
resolutive-regular point,
weak Kellogg property.
}


\section{Introduction}

In this paper we study boundary regularity for \p-harmonic
functions on metric spaces, including $\R^n$ as an important case.
Such studies have earlier been done with respect to the given metric boundary,
but the novelty here is that we consider boundary regularity
with respect to the Mazurkiewicz boundary and other compactifications.

This builds on the earlier work for the Dirichlet problem
with respect to the Mazurkiewicz boundary, by 
Bj\"orn--Bj\"orn--Shanmugalingam~\cite{BBSdir},
and more recently with respect to arbitrary compactifications,
by Bj\"orn--Bj\"orn--Sj\"odin~\cite{BBSjodin}.
Boundary regularity was however not considered therein.

To be more precise, let $1 < p < \infty$ and let
$X$ be a complete metric space 
equipped with a doubling measure $\mu$ supporting a \p-Poincar\'e inequality.
$X$ can e.g.\ be unweighted $\R^n$,
in which case a function $u$ is \emph{\p-harmonic}
 if it is a continuous weak solution
of the \p-Laplace equation
\[
    \Div(|\nabla u|^{p-2} \nabla u)=0.
\]
The definition of \p-harmonic
functions on metric spaces is more involved, see Section~\ref{sect-pharm}.

Let $\Om$ be a bounded domain.
(If $X$ is bounded we also require that
$\Cp(X \setm \Om)>0$.)
The \emph{Mazurkiewicz distance} $\dM$ on $\Om$ is defined by
\[
     \dM(x,y) =\inf_E \diam E,
\]
where the infimum is taken over all connected sets $E \subset \Om$
containing $x,y \in \Om$.
(The Mazurkiewicz distance was first used by
Mazurkiewicz~\cite{mazurkiewicz16} in 1916,
but goes under different names in the literature,
see Remark~4.2 in \cite{BBSdir}.)
The \emph{Mazurkiewicz boundary} $\bdym \Om$ is 
the boundary of $\Om$ in the completion
of $(\Om,\dM)$.
For instance, in the slit disc 
$B(0,1) \setm [0,1] \subset \R^2$
this gives two  boundary points corresponding to each point in the slit
(but for the tip),
while for smooth domains $\bdym \Om = \bdy \Om$.

Assume that the completion of $(\Om,\dM)$
is compact, which happens if and only if $\Om$ is 
finitely connected at the boundary, see Section~\ref{sect-compactification}.
Then, to each point in the given metric boundary $\bdy \Om$
there corresponds one or more (at most countably many)
points in the Mazurkiewicz boundary $\bdym \Om$,
while conversely to every point in $\bdym \Om$ there
corresponds a unique point in $\bdy \Om$,
see Bj\"orn--Bj\"orn--Shanmugalingam~\cite{BBSmbdy}.
There is therefore a natural projection $\Phi: \bdym \Om \to \bdy \Om$
between these boundaries.

The following is the first new result, and the key tool to obtaining
the Kellogg property.

\begin{thm} \label{thm-reg-bdym}
Assume that $\Om$ is a bounded domain which is
finitely connected at the boundary.
Let $x_0 \in \bdy \Om$.
Then the following are equivalent\/\textup{:}
\begin{enumerate}
\item 
The point $x_0$ is an irregular boundary point with respect to $\Om$.
\item 
There is at least one irregular boundary point in $\Phi^{-1}(x_0)$
with respect to $\Om^M$.
\item 
There is exactly  one irregular boundary point in $\Phi^{-1}(x_0)$
with respect to $\Omm$.
\end{enumerate}
\end{thm}

This leads directly to the following consequence.

\begin{thm}  \label{thm-Kellogg-bdym}
\textup{(The Kellogg property)}
Assume that $\Om$ is a bounded domain which is
finitely connected at the boundary.
Let $\Irr^M$ be the set of irregular boundary points
with respect to $\Omm$.
Then the capacity $\bCp(\Irr^M,\Omm)=0$.
\end{thm}

Note that by Proposition~\ref{prop-cap-all-bdy},
$\bCp(\bdym\Om,\Omm)>0$,
so the Kellogg property
is never trivial.
The following uniqueness result is also new.

\begin{thm} \label{thm-unique-bdym}
Assume that $\Om$ is a bounded domain which is
finitely connected at the boundary.
Let $f \in C(\bdym \Om)$. Then there exists
a unique bounded \p-harmonic function $u$ on $\Om$ such that
\[
    \lim_{\Om \ni y \dMto x} u(y) = f(x) \quad 
\text{for $\bCp(\,\cdot\,;\Omm)$-q.e. } x \in \bdym \Om.
\]
Moreover, $u$ equals  the Perron solution $\Hpind{\Omm} f$.
\end{thm}

We are also able to show that  boundary regularity
is a local property for the Mazurkiewicz boundary
in the following sense.

\begin{thm} \label{thm-local-prop-bdym}
Assume that $\Om$ is a bounded domain which is
finitely connected at the boundary.
Let $\xh \in \bdym \Om$ and let 
$G $ be an $\clOmm$-neighbourhood of $\xh$.

Then $\xh$ is regular with respect to $\Omm$
if and only if it is regular with respect to  $G^{\Omm}$,
where $G^{\Omm}$ is $G$ equipped with the boundary inherited from $\clOmm$.
\end{thm}

Throughout the paper we also study to what extent such results are true
for other compactifications of $\Om$. The
details describing the different results and cases are quite
involved.

Boundary regularity for \p-harmonic functions with respect to the given metric 
boundary
has been studied for a long period, especially on $\R^n$.
The first significant result was Maz$'$ya's~\cite{Maz70}
sufficiency part of the Wiener criterion in 1970.
Later on the full Wiener criterion was obtained
in various situations including weighted $\R^n$ and for Cheeger \p-harmonic
functions on metric spaces,
see \cite{HeKiMa}, \cite{KiMa94}, \cite{Lind-Mar}, \cite{Mikkonen}
and \cite{JB-Matsue}.
The full Wiener criterion remains open (for the given metric boundary)
in the generality considered here, but the sufficiency has been obtained, 
see \cite{BMS} and \cite{JB-pfine},
and a weaker necessity condition, see \cite{JBCalcVar}.

In the nonlinear potential theory, the Kellogg property  was first obtained by
Hedberg~\cite{Hedb} and  Hedberg--Wolff~\cite{HedWol} on $\R^n$
(see also Kilpel\"ainen~\cite{Kilp89}),
who also obtained the more general fine Kellogg property and the even more
general Choquet property.
The Kellogg property was
extended to homogeneous spaces by Vodop{\cprime}\-yanov~\cite{Vodopyanov89}
(who also obtained the fine Kellogg property),
to weighted $\R^n$ by
Hei\-no\-nen--Kil\-pe\-l\"ai\-nen--Martio~\cite{HeKiMa},
to subelliptic equations by  
Markina--Vodop{\cprime}yanov~\cite{MV2},
and to metric spaces by 
Bj\"orn--Bj\"orn--Shan\-mu\-ga\-lin\-gam~\cite{BBS}.
The fine Kellogg property, as well as the Choquet property,
was deduced on metric spaces by Bj\"orn--Bj\"orn--Latvala~\cite{BBL2}.
Other aspects of boundary regularity and boundary behaviour
for \p-harmonic functions
have been studied in
\cite{ABbarrier}--\cite{BB2},
\cite{GLM86}, 
\cite{HeKiMa}, 
\cite{Kilp89},
\cite{KiLi00}
and 
\cite{MV2}.

The Wiener criterion characterizes the regularity of a boundary point
using the complement of the domain (beyond the boundary).
Many other results, such as the Kellogg property, do not directly involve the
complement, but the proofs of most boundary regularity results do use
the complement (beyond the boundary) in significant ways.

In our situation we have a boundary of the domain,
but no complement beyond that.
Thus most of the techniques used to study boundary regularity
with respect to the given metric boundary are not available to us.
Instead we will mainly depend on comparing boundary regularity
between different boundaries.
In particular, most of our stronger results are for boundaries
larger than the given metric boundary.

We also give several counterexamples showing that the specific assumptions
in our results are at least to some extent essential.
These include two examples where the Kellogg property fails.

\begin{ack}
The author was supported by the Swedish Research Council,
grants 621-2007-6187, 621-2011-3139  and 2016-03424.
The idea to study resolutive-regularity is due to 
Tomas Sj\"odin (private communication).
\end{ack}

\section{Notation and preliminaries}
\label{sect-prelim}

We will need quite a bit of notation, which we will introduce
in this and the next two sections. 
We will be brief, see
Bj\"orn--Bj\"orn--Shanmugalingam~\cite{BBSdir} 
and Bj\"orn--Bj\"orn--Sj\"odin~\cite{BBSjodin}
for more details.
Proofs of the results in 
this section can be found in the monographs
Bj\"orn--Bj\"orn~\cite{BBbook} and
Heinonen--Koskela--Shanmugalingam--Tyson~\cite{HKSTbook}.

We assume throughout the  paper that $1 < p<\infty$ 
and that $X=(X,d,\mu)$ is a metric space equipped
with a metric $d$ and a positive complete  Borel  measure $\mu$ 
such that $0<\mu(B)<\infty$ for all balls $B \subset X$.

We will only consider curves which are nonconstant, compact
and 
rectifiable (i.e.\ have finite length), and thus each curve can 
be parameterized by its arc length $ds$. 
A property is said to hold for \emph{\p-almost every curve}
if it fails only for a curve family $\Ga$ with zero \p-modulus, 
i.e.\ there exists $0\le\rho\in L^p(X)$ such that 
$\int_\ga \rho\,ds=\infty$ for every curve $\ga\in\Ga$.

Following Koskela--MacManus~\cite{KoMc} (see also 
Heinonen--Koskela~\cite{HeKo98}) we
introduce weak upper gradients as follows.

\begin{deff} \label{deff-ug}
A measurable function $g : X \to [0,\infty]$  is a \emph{\p-weak upper gradient} 
of a function $f: X \to \eR:=[-\infty,\infty]$
if for \p-almost all  curves  
$\gamma : [0,l_{\gamma}] \to X$,
\begin{equation} \label{ug-cond}
        |f(\gamma(0)) - f(\gamma(l_{\gamma}))| \le \int_{\gamma} g\,ds,
\end{equation}
where the left-hand side is considered to be $\infty$ 
whenever at least one of the 
terms therein is infinite.
\end{deff}

If $f$ has a \p-weak upper gradient in $\Lploc(X)$, then
it has an a.e.\ unique \emph{minimal \p-weak upper gradient} $g_f \in \Lploc(X)$
in the sense that for every \p-weak upper gradient $g \in \Lploc(X)$ of $f$ we have
$g_f \le g$ a.e., see Shanmugalingam~\cite{Sh-harm}.
Following Shanmugalingam~\cite{Sh-rev}, 
we define a version of Sobolev spaces on the metric space $X$.

\begin{deff} \label{deff-Np}
For a measurable function $f: X\to \eR$, let 
\[
        \|f\|_{\Np(X)} = \biggl( \int_X |f|^p \, d\mu 
                + \inf_g  \int_X g^p \, d\mu \biggr)^{1/p},
\]
where the infimum is taken over all \p-weak upper gradients $g$ of $f$.
The \emph{Newtonian space} on $X$ is 
\[
        \Np (X) = \{f: \|f\|_{\Np(X)} <\infty \}.
\]
\end{deff}
\medskip
The space $\Np(X)/{\sim}$, where  $f \sim h$ if and only if $\|f-h\|_{\Np(X)}=0$,
is a Banach space,  see \cite{Sh-rev}.
We also define
\[
   \Dp(X)=\{f : f \text{ is measurable and  has a \p-weak  upper gradient
     in }   L^p(X)\}.
\]
In this paper we assume that functions in $\Np(X)$
and $\Dp(X)$
 are defined everywhere (with values in $\eR$),
not just up to an equivalence class in the corresponding function space.
For a measurable set $E\subset X$, the Newtonian space $\Np(E)$ is defined by
considering $(E,d|_E,\mu|_E)$ as a metric space in its own right.
We say  that $f \in \Nploc(E)$ if
for every $x \in E$ there exists a ball $B_x\ni x$ such that
$f \in \Np(B_x \cap E)$.

\begin{deff} \label{deff-sobcap}
The (Sobolev) \emph{capacity} of an arbitrary set $E\subset X$ is
\[
\Cp(E) = \inf_u\|u\|_{\Np(X)}^p,
\]
where the infimum is taken over all $u \in \Np(X)$ such that
$u\geq 1$ on $E$.
\end{deff}

The measure  $\mu$  is \emph{doubling} if 
there exists a \emph{doubling constant} $C>0$ such that 
\begin{equation*}
        0 < \mu(2B) \le C \mu(B) < \infty
\end{equation*}
for all balls 
$B=B(x_0,r):=\{x\in X: d(x,x_0)<r\}$ in~$X$,
where $\lambda B=B(x_0,\lambda r)$.  

\begin{deff} \label{def.PI.}
$X$ supports a \emph{\p-Poincar\'e inequality} if
there exist constants $C>0$ and $\lambda \ge 1$
such that for all balls $B \subset X$, 
all integrable functions $f$ on $X$ and all 
\p-weak upper gradients $g$ of $f$, 
\[ 
        \vint_{B} |f-f_B| \,d\mu 
        \le C \diam(B) \biggl( \vint_{\lambda B} g^{p} \,d\mu \biggr)^{1/p},
\] 
where $ f_B 
 :=\vint_B f \,d\mu 
:= \int_B f\, d\mu/\mu(B)$.
\end{deff}

In this paper neighbourhoods are always open
and
continuous functions are real-valued,
whereas semicontinuous functions may take the values $\pm \infty$.

\section{Compactifications and the capacity
\texorpdfstring{$\bCp$}{}}
\label{sect-compactification}

\begin{deff}
Let $\Om$ be a locally compact noncompact metric space.
A couple $(\partial' \Omega,\tau)$
is said to \emph{compactify} $\Omega$ if
$\partial' \Omega$ is a set with $\partial' \Omega \cap \Omega = \emptyset$
and $\tau$ is a Hausdorff topology on 
$\clOmprim:=\Omega \cup \partial' \Omega$ such that
\begin{enumerate}
\renewcommand{\theenumi}{\textup{(\roman{enumi})}}%
\item $\clOmprim$  is compact with respect to $ \tau$;
\item \label{b} 
$\Omega$ is dense in $\clOmprim$ with respect to $\tau$;
\item \label{c}
the topology induced on $\Omega$ by $\tau$  coincides with
the given topology on $\Om$.
\end{enumerate}
The space $\clOmprim$ with the topology $\tau$ is 
a \emph{compactification} of $\Omega$.
\end{deff}

Since $\clOmprim$ is a compact Hausdorff space it is normal.

\medskip

\emph{We assume from now on that $X$ is a complete
metric space supporting a \p-Poincar\'e inequality,
that $\mu$ is doubling, and that\/ $1<p<\infty$.
We also assume that\/ $\Om$ is a 
nonempty bounded open set
such that $\Cp(X \setm \Om)>0$,
and that 
$\clOmj=\Om^j \cup \bdyj \Om$, $j=1,2$, are compactifications of 
$\Om$, where $\Om^j= \Om$ with the intended boundary $\bdyj \Om$
and where the topologies on $\clOmj$ are denoted by $\tauj$.
Furthermore, we reserve $\partial \Omega$ and $\clOm$ 
for the given metric boundary 
and closure induced by $X$ on $\Omega$.}

\medskip

As $\mu$  is doubling and $X$ is complete, it follows
that $X$ is proper (i.e.\ all closed bounded sets are compact).

We define 
$
\partial^1 \Omega \prec \partial^2 \Omega
$
to mean that there is a continuous mapping, 
which is called \emph{projection},
\[
 \Phi: \clOmtwo \longrightarrow \clOmone
\quad \text{with } \Phi|_\Om=\id.
\]

An example of a compactification is the Mazurkiewicz completion
discussed in the introduction,
which is a compactification of $\Om$ if and only if 
$\Om$ is a domain which is finitely
connected at the boundary (in the following sense), by Theorem~1.1 in
Bj\"orn--Bj\"orn--Shanmugalingam~\cite{BBSmbdy} or
Theorem~1.3.8 in Karmazin~\cite{Ka}.

\begin{deff} \label{deff-fin-conn}
A bounded domain $\Om \subset X$ is \emph{finitely connected at the boundary} 
if for every  $x\in \bdy \Om$ and  $r>0$ there is an open set $G$ (in $X$) 
such that $x \in G \subset B(x,r)$
and $G \cap \Om$ has only finitely many components.
\end{deff}

In addition to the Sobolev capacity mentioned above, we will also need
the following capacity. 
It was introduced for $\clOm$  and $\clOmm$ by
Bj\"orn--Bj\"orn--Shan\-mu\-ga\-lin\-gam~\cite{BBSdir},
and generalized to arbitrary compactifications as here
in 
Bj\"orn--Bj\"orn--Sj\"odin~\cite[Definition~4.1]{BBSjodin}.
A similar capacity was considered in 
Kilpel\"ainen--Mal\'y~\cite{KilMaGenDir}.

\begin{deff} \label{deff-bCp}
For $E \subset \clOmone$  let
\[
     \bCp(E;\Omone)= \inf_{u \in \A_E(\Omone)} \|u\|_{\Np(\Om)}^p,
\]
where $u \in \A_E(\Omone)$ if $u \in \Np(\Om)$ is such that
\begin{equation} \label{eq-deff-bCp-1}
  u \ge 1 
  \quad  \text{on } E \cap \Om
\end{equation}
and 
\begin{equation} \label{eq-deff-bCp}
      \liminf_{\Om \ni y \toone x} u(y) \ge 1
	\quad \text{for all } x \in E \cap \bdyone \Om.
\end{equation}
\end{deff}

When proving the Kellogg property we will need the following lemma.

\begin{lem} \label{lem-bCp-Phi}
Assume that $\bdyone \Om \prec \bdytwo \Om$
and let $\Phi: \clOmtwo \to \clOmone$ be the projection.
Let $E \subset \bdyone \Om$.
Then 
$\bCp(E;\Omone) = \bCp(\Phi^{-1}(E),\Omtwo)$.
\end{lem}

We will use nets to study convergence in our compactifications,
see e.g.\ Pedersen~\cite{pedersen} for the key results on nets.

\begin{proof}
We first show that $\A_E \subset \A_{\Phi^{-1}(E)}$.
Take $u \in \A_E$.
Let $\xh \in \Phi^{-1}(E)\cap\bdytwo\Om$ and take a 
net $y_\la \in \Om$ such that $y_\la \totwo \xh$. 
Since $\Phi$ is continuous it follows that
$y_\la \toone \Phi(\xh)$.
Hence, because $\Phi(\xh) \in E\cap\bdyone\Om$,
we see that 
$     \liminf_{\la} u(y_\la) \ge 1$, 
and thus $u \in \A_{\Phi^{-1}(E)}$.

Conversely, assume that $u \in \A_{\Phi^{-1}(E)}$.
Let $x \in E\cap\bdyone\Om$.
Assume that there is a 
net $y_\la \in \Om$ such that $y_\la \toone x$
and $\liminf_\la u(y_\la) < 1$.
By taking a subnet we may assume that $\lim_\la u(y_\la)$ exists
and is less than $1$.
Then there is a further subnet
$y_{\mu}$ which  converges to some point $\xh \in\Omtwo$.
Hence $\lim_{\mu} u(y_{\mu}) < 1$.
As $\Phi$ is continuous we see that $y_{\mu} \toone \Phi(\xh)$,
and thus $\Phi(\xh)=x$.
But together with $\lim_{\mu} u(y_{\mu}) < 1$,
this contradicts the assumption $u \in \A_{\Phi^{-1}(E)}$.
Hence 
\[
\liminf_{y \toone x} u(y) \ge 1
\]
 and 
$u \in A_E$.

Therefore the infima defining
the two capacities are taken over the same set and 
the two capacities agree.
\end{proof}

As a consequence we can obtain the following result,
which shows that the Kellogg property is not seeking something trivial.

\begin{prop} \label{prop-cap-all-bdy}
$\bCp(\bdyone\Om;\Omone)>0$.
\end{prop}

\begin{proof}
Let $\clOmtwo$ be the one-point compactification of $\Om$.
We will first show that $\bCp(\bdy\Omtwo;\Omtwo)>0$.
Assume on the contrary that $\bCp(\bdy\Omtwo;\Omtwo)=0$.
Then for each $j=1,2,\ldots$, there is $u_j \in \A_{\bdy\Omtwo}(\Omtwo)$
with $0 \le u_j \le 1$ and $\|u_j\|_{\Np(\Om)} < 1/j$.
Let
\[
v_j =\begin{cases}
    (1-2u_j)_\limplus & \text{in } \Om, \\
    0 & \text{in } X \setm \Om,
   \end{cases}
   \quad \text{and} \quad
   v(x)=\limsup_{j \to \infty} v_j (x), \ x \in X.
\]
By \eqref{eq-deff-bCp}, the set 
$(X \setm \Om) \cup \bigl\{x \in \Om : u_j(x) > \tfrac12\bigr\}$
contains a neighbourhood of
$\bdy \Om$, and therein $v_j \equiv 0$.
As $v_j \in \Np(\Om)$ 
we conclude that  $v_j \in \Np(X)$.
Moreover, $\{v_j\}_{j=1}^\infty$ is a Cauchy sequence in $\Np(X)$
and thus, by Corollary~1.72 in \cite{BBbook}, it has a subsequence
which converges q.e.\ to $v \in \Np(X)$.
We get directly that $v \equiv 0$ in $X \setm \Om$.
Moreover, if $0 <\de <\tfrac{1}{2}$, then
\[
  \mu(\{x \in \Om : v_j(x)  \le 1 -2 \de\}
   = \mu(\{x \in \Om : u_j(x)  \ge  \de\}
  \le \frac{1}{\de^p} \int_\Om u_j^p \, d\mu 
  \le \frac{1}{(j\de)^p},
\]
which tends to $0$ as $j \to \infty$.
It follows that $v=1$ a.e.\ in $\Om$.

We next need to consider two cases separately.
First, if $\mu(X \setm \Om)=0$, then
$v$ equals $1$ a.e.\  but not q.e.\ in $X$ (because $\Cp(X \setm \Om)>0$), 
contradicting
Proposition~1.59 in \cite{BBbook}.

On the other hand, if $\mu(X \setm \Om)>0$, then
we let $B$ be a large enough ball so that $\Om \subset B$ and
$\mu(B \setm \Om)>0$.
By Corollary~2.21 in \cite{BBbook}, we know that $g_v = 0$ a.e.
Hence, by the \p-Poincar\'e inequality,
\[
    0 <         \vint_{B} |v-v_B| \,d\mu 
    \le C \diam(B) \biggl( \vint_{\lambda B} g_v^{p} \,d\mu \biggr)^{1/p}
    = 0,
\]
a contradiction.

Thus we must have $\bCp(\bdy\Omtwo;\Omtwo)>0$,
and by 
Lemma~\ref{lem-bCp-Phi},
\[
    \bCp(\bdyone\Om;\Omone)= \bCp(\bdytwo\Om;\Omtwo)>0.
    \qedhere
\]
\end{proof}

\section{\texorpdfstring{\p}{p}-harmonic functions and Perron solutions}
\label{sect-pharm}

A function $u \in \Nploc(\Om)$ is a
\emph{\textup{(}super\/\textup{)}minimizer} in $\Om$
if 
\[ 
      \int_{\phi \ne 0} g^p_u \, d\mu
           \le \int_{\phi \ne 0} g_{u+\phi}^p \, d\mu
           \quad \text{for all (nonnegative) } \phi \in \Np_0(\Om).
\] 
A \emph{\p-harmonic function} is a continuous minimizer.

A function $u : \Om \to (-\infty,\infty]$ is \emph{superharmonic}
in $\Om$ if 
\begin{enumerate}
\renewcommand{\theenumi}{\textup{(\roman{enumi})}}%
\renewcommand{\labelenumi}{\theenumi}%
\item
$u$ is not identically $\infty$ in any component of $\Om$; 
\item
$\min\{u,k\}$ is an lsc-regularized superminimizer for all $k \in \R$,
\end{enumerate}
where $v$ is 
\emph{lsc-regularized} if
\[ 
 v(x):=\essliminf_{y\to x} v(y):= \lim_{r \to 0} \essinf_{B(x,r)} v.
\] 

By Theorem~6.1 in Bj\"orn~\cite{ABsuper}
(or \cite[Theorems~9.24 and 14.10]{BBbook}),
this definition of superharmonicity is equivalent to  the ones used
both in the Euclidean and metric space literature, e.g.\
in
Hei\-no\-nen--Kil\-pe\-l\"ai\-nen--Martio~\cite{HeKiMa}, 
Kinnunen--Martio~\cite{KiMa02} and Bj\"orn--Bj\"orn~\cite{BBbook}. 

We are now ready to introduce the Perron solutions
with respect to $\Omone$.
We follow Bj\"orn--Bj\"orn--Sj\"odin~\cite{BBSjodin},
although therein Perron solutions were only defined in domains.

\begin{deff}   \label{def-Perron}
Given 
$f : \bdyone \Om \to \eR$, let $\UU_f(\Omone)$ be the set of all 
superharmonic functions $u$ on $\Om$, bounded from below,  such that 
\begin{equation} \label{eq-def-Perron} 
	\liminf_{\Om \ni y \toone x} u(y) \ge f(x) 
\end{equation} 
for all $x \in \bdyone \Om$.
The \emph{upper Perron solution} of $f$ is then defined to be
\[ 
    \uHpind{\Omone} f (x) = \inf_{u \in \UU_f(\Omone)}  u(x), \quad x \in \Om,
\]
while the \emph{lower Perron solution} of $f$ is defined by
\[ 
    \lHpind{\Omone} f = - \uHpind{\Omone} (-f).
\]
If $\uHpind{\Omone} f = \lHpind{\Omone} f$
and it is real-valued, then we let $\Hpind{\Omone} f := \uHpind{\Omone} f$
and $f$ is said to be \emph{resolutive} with respect to $\Omone$.

Furthermore, let $\DU_f(\Omone)=\UU_f(\Omone) \cap \Dp(\Om)$,
and define the \emph{Sobolev--Perron solutions} of $f$ by
\begin{equation} \label{eq-Sp} 
    \uSpind{\Omone} f (x) = \inf_{u \in \DU_f(\Omone)}  u(x),\ x \in \Om,
\quad    \text{and}  \quad
    \lSpind{\Omone} f  = -\uSpind{\Omone} (-f).
\end{equation}
If $\uSpind{\Omone} f = \lSpind{\Omone} f$ and it is real-valued, 
then 
$f$ is said to be \emph{Sobolev-resolutive} with respect to $\Omone$.
The boundary $\bdyone \Om$ is \emph{\textup{(}Sobolev\/\textup{)}-resolutive} if all 
functions $f \in C (\bdyone \Om)$ are (Sobolev)-resolutive.
\end{deff}

In every component of $\Om$ the upper/lower (Sobolev)--Perron solutions
are \p-harmonic or identically $\pm \infty$,
see Theorem~4.1 in Bj\"orn--Bj\"orn--Shan\-mu\-ga\-lin\-gam~\cite{BBS2}
(or Theorem~10.10 in \cite{BBbook}). The Sobolev--Perron
solutions were introduced in Bj\"orn--Bj\"orn--Sj\"odin~\cite{BBSjodin};
we will only use them in Corollary~\ref{cor-Sobolev-reg}.

The given metric boundary $\bdy \Om$ is resolutive by
Theorem~6.1 in \cite{BBS2}
(or Theorem~10.22 in \cite{BBbook}).
It is also Sobolev-resolutive by
Theorem~6.4 and Proposition~7.3 in \cite{BBSjodin}.
If $\Om$ is finitely connected at the boundary,
then $\bdym \Om$ is resolutive by
Theorem~8.2 in Bj\"orn--Bj\"orn--Shanmugalingam~\cite{BBSdir}.
Also $\bdym \Om$ is Sobolev-resolutive (if
$\Om$ is finitely connected at the boundary), which again
follows from Theorem~6.4 and Proposition~7.3 in \cite{BBSjodin}
since continuous functions on $\bdym \Om$ can be uniformly
approximated by Lipschitz functions on $\clOmm$.

The following results from 
\cite{BBSjodin}
will be important for us.

\begin{prop} \label{prop-lHp<=uHp}
\textup{(\cite[Corollary~6.3]{BBSjodin})}
If $f:\bdyone\Om\to\eR$, then  
$  \lHpind{\Omone} f      \le \uHpind{\Omone} f$.
\end{prop}

\begin{thm} \label{thm-res-prec}
\textup{(\cite[Theorem~6.7 and Proposition~8.2]{BBSjodin})}
Assume that $\partial^1 \Omega \prec \partial^2 \Omega$
and let $\Phi : \clOmtwo \to \clOmone$ 
denote the projection.
If $f : \partial^1 \Omega \rightarrow \eR$, then  
\begin{equation} \label{eq-B}
\uP_{\Om^1}f = \uP_{\Om^2}({f \circ \Phi}) 
\quad \text{and} \quad
\lP_{\Om^1}f = \lP_{\Om^2}({f \circ \Phi}).
\end{equation}
In particular, if $\partial^2 \Omega$ is
resolutive then so is $\partial^1 \Omega$.
\end{thm}

One consequence of Theorem~\ref{thm-res-prec}
is that (almost always) there are plenty
of resolutive functions.

\section{Boundary regularity}
\label{sect-bdy-reg}

Resolutivity will play an important role in several of our boundary
regularity results. One possibility would have been to restrict our
attention to resolutive boundaries. Here we have instead
chosen a more general approach introducing both regular and resolutive-regular
boundary points. 
The idea of studying resolutive-regularity is due to
Sj\"odin (private communication).

\begin{deff} \label{deff-reg-pt}
A boundary point $x_0 \in \bdyone \Om$ is 
\emph{\textup{(}resolutive\/\textup{)}-regular}
if
\[
              \lim_{\Om \ni y \toone x_0} \uHpind{\Omone} f(y) = f(x_0)
\]
for all (resolutive)
$f\in C(\bdyone \Om)$.
Otherwise, it is \emph{\textup{(}resolutive\/\textup{)}-irregular}.

We also say that $\partial^1 \Omega$ is
\emph{\textup{(}resolutive\/\textup{)}-regular}
if all its boundary points are (re\-so\-lu\-tive)-regular.
\end{deff}

Note that if all continuous functions are resolutive, then 
regularity and resolutive-regularity
of course coincide.
Example~\ref{ex-4-p>2} shows that this is not true in general.
The following result shows that we can equivalently replace
$\lim$ by $\limsup$ and $=$ by $\le$ in Definition~\ref{deff-reg-pt}.

\begin{lem} \label{lem-res-reg}
A boundary point $x_0 \in \bdyone \Om$ is\/ 
\textup{(}resolutive\/\textup{)}-regular
if and only if 
\begin{equation} \label{eq-res-reg}
              \limsup_{\Om \ni y \toone x_0} \uHpind{\Omone} f(y) \le f(x_0)
\end{equation}
for all \textup{(}resolutive\/\textup{)} 
$f\in C(\bdyone \Om)$.
\end{lem}

\begin{proof}
Assume that \eqref{eq-res-reg} holds for all (resolutive) $f\in C(\bdyone \Om)$.
Let $f:\bdyone \Om \to \R$ be continuous (and resolutive).
Then also $-f$ is continuous (and resolutive).
Hence, by Proposition~\ref{prop-lHp<=uHp},
\[
    f(x_0) \le 
\liminf_{\Om \ni y \toone x_0} \lHpind{\Omone} f(y) 
\le \liminf_{\Om \ni y \toone x_0} \uHpind{\Omone} f(y) 
\le \limsup_{\Om \ni y \toone x_0} \uHpind{\Omone} f(y) 
   \le f(x_0),
\]
and therefore
\[
   \lim_{\Om \ni y \toone x_0} \uHpind{\Omone} f(y) 
   = f(x_0).
\]
Thus $x_0$ is 
(resolutive)-regular.
The converse is trivial.
\end{proof}

The following result shows that the boundary
regularity classification (into regular and irregular boundary points) can be useful
also for noncontinuous boundary data.
As regularity and resolutive-regularity are different, by
Example~\ref{ex-4-p>2}, the latter cannot be characterized
in a similar fashion.

\begin{prop} \label{prop-reg-1}
Let $x_0 \in \bdyone \Om$.
Then the following are equivalent\/\textup{:}
\begin{enumerate}
\renewcommand{\theenumi}{\textup{(\roman{enumi})}}%
\item \label{reg-1} 
The point $x_0$ is a regular boundary point.
\item \label{cont-x0-perron-1}
It is true that 
\[
        \lim_{\Om \ni y \toone x_0}  \uHpind{\Omone} f(y) = f(x_0)
\]
for all bounded $f :\bdyone\Om \to \R$
which are continuous at $x_0$.
\item \label{semicont-x0-perron-1}
It is true that 
\begin{equation} \label{eq-usc}
        \limsup_{\Om \ni y \toone x_0}  \uHpind{\Omone} f(y) \le f(x_0)
\end{equation}
for all functions $f :\bdyone\Om \to [-\infty,\infty)$ which are
bounded from above on $\bdyone \Om$ and 
upper semicontinuous at $x_0$.
\end{enumerate}
\end{prop}

\begin{proof}
\ref{cont-x0-perron-1} $\imp$ \ref{reg-1}
This is trivial.

\ref{reg-1} $\imp$ \ref{semicont-x0-perron-1} 
Let $A >f(x_0)$ be real and $M=\sup_{\bdyone\Om} f_\limplus$.
Since $f$ is upper semicontinuous at $x_0$ we can find a 
$\tauone$-neighbourhood $G \subset \bdyone \Om$ containing $x_0$ such that
$f <A$ in $G$.
By Tietze's extension theorem,
we can also find $h \in C(\bdyone \Om)$ 
such that $h(x_0)=A$, $h=M$ on $\bdyone \Om \setm G$
and $h \ge A$ everywhere.
Then $f \le h$ on $\bdyone \Om$,
and thus
\[
     \limsup_{\Om \ni y \toone x_0} \uHpind{\Omone} f(y) 
      \le \lim_{\Om \ni y \toone x_0} \uHpind{\Omone} h(y) 
      = h(x_0)=A,
\]
from which \eqref{eq-usc} follows after taking infiumum
over all $A > f(x_0)$.

\ref{semicont-x0-perron-1} $\imp$ \ref{cont-x0-perron-1} 
Applying \ref{semicont-x0-perron-1} to $-f$ yields
\begin{align*}
     \liminf_{\Om \ni y \toone x_0} \uHpind{\Omone} f(y) 
     & \ge \liminf_{\Om \ni y \toone x_0} \lHpind{\Omone} f(y) \\
     & = - \limsup_{\Om \ni y \toone x_0} \uHpind{\Omone}(-f)(y) 
      \ge - (-f(x_0)) = f(x_0).
\end{align*}
Together with \eqref{eq-usc} this gives the 
desired conclusion.
\end{proof}

\begin{prop} \label{prop-irreg-1=>2}
Assume that $\bdyone \Om \prec \bdytwo \Om$
and let $\Phi: \clOmtwo \to \clOmone$ be the projection.
If  $x_0 \in \bdyone \Om$ is\/ \textup{(}resolutive\/\textup{)}-irregular
with respect to $\Omone$,
then there is at least one\/ \textup{(}resolutive\/\textup{)}-irregular
boundary point in $\Phi^{-1}(x_0)$
with respect to $\Om^2$.
\end{prop}

\begin{proof}
Since $x_0$ is \textup{(}resolutive\/\textup{)}-irregular 
there is, due to Lemma~\ref{lem-res-reg},
a (resolutive) function  $f \in C(\bdy^1 \Om)$ such that
\[
             M:= \limsup_{\Om \ni y \toone x_0} \uHpind{\Omone} f(y) > f(x_0).
\]
We may assume that $M=2$ and $f(x_0)=0$.
Let $G=\{y \in \Om : \uHpind{\Omone} f(y) > 1\}$, which is a nonempty open subset of
$\Om$ such that $x_0 \in \clGone$, where $\clGone$ is the
closure of $G$ within $\clOmone$.
In particular there is a net $\{x_\la\}_\la$ such that 
$G \ni x_\la \toone x_0$.
As $G$ is not closed in $\clOmone$ it is not compact.

By compactness of $\clOmtwo$, there is a $\tautwo$-convergent subnet $\{x_\mu\}_\mu$,
with $\tautwo$-limit $\xh$.
It follows that $\Phi(\xh)=x_0$.
Set $h=f \circ \Phi$.
By Theorem~\ref{thm-res-prec},
$\uHpind{\Om^2} h =\uHpind{\Om^1} f $
(and $h$ is resolutive if $f$ is).
Thus
\[
\limsup_{\Om \ni y \totwo \xh} \uHpind{\Omtwo} h(y)
= \limsup_{\Om \ni y \totwo \xh} \uHpind{\Omone} f(y)
\ge  \limsup_{\mu} \uHpind{\Omone} f(x_\mu)
\ge 1.
\]
As $h(\xh)=f(x_0)=0$ and $h$ is continuous (and resolutive if $f$ is), this shows that
$\xh$ is (resolutive)-irregular with respect to $\Omtwo$.
\end{proof}

\section{The proof of Theorem~\ref{thm-reg-bdym}}

We are now ready to consider our generalization of
Theorem~\ref{thm-reg-bdym}.
For this we need an additional assumption, which we now define.

\begin{deff} \label{deff-nicely}
Assume that $\bdyone \Om \prec \bdytwo \Om$
and let $\Phi: \clOmtwo \to \clOmone$ be the projection.
We say that  $x_0 \in \bdyone \Om$ \emph{splits nicely} if
every $\xh \in \Phi^{-1}(x_0)$ has arbitrarily small
$\tautwo$-neighbourhoods $U$ such that 
$\bdytwo U \cap \Phi^{-1}(x_0)=\emptyset$.
\end{deff}

Note that if $\Phi^{-1}(x_0)$ is finite,
or if $\clOmtwo$ is metrizable and $\Phi^{-1}(x_0)$ is at most
countable, then $x_0$ splits nicely.
In particular, all points in $\bdy \Om$ split nicely with
respect to the Mazurkiewicz boundary $\bdym \Om$.
When $\clOmtwo$ is nonmetrizable, we do not know if
$x_0$ always splits nicely whenever
$\Phi^{-1}(x_0)$ is countable.

\begin{thm} \label{thm-reg-Phi}
Assume that $\bdy \Om \prec \bdytwo \Om$,
where $\bdy \Om$ is the given metric boundary.
Let $\Phi: \clOmtwo \to \clOm$ be the projection.
Assume that $x_0 \in \bdy \Om$ splits nicely.
Then the following are equivalent\/\textup{:}
\begin{enumerate}
\item \label{i-Om-reg}
The point $x_0$ is irregular  with respect to $\Om$.
\item \label{i-Om}
The point $x_0$ is resolutive-irregular with respect to $\Om$.
\renewcommand{\theenumi}{\textup{(\alph{enumi})}}%
\item \label{i-Omtwo-at-least}
There is at least one resolutive-irregular boundary point in $\Phi^{-1}(x_0)$
with respect to $\Omtwo$.
\item \label{i-Omtwo-exactly}
There is exactly  one resolutive-irregular boundary point in $\Phi^{-1}(x_0)$
with respect to $\Omtwo$.
\end{enumerate}
\end{thm}

On the way to proving Theorem~\ref{thm-reg-Phi}
we first obtain the following result,
which generalizes both \ref{i-Omtwo-at-least} $\imp$ \ref{i-Om}
and \ref{i-Omtwo-exactly} $\imp$ \ref{i-Om} in  Theorem~\ref{thm-reg-Phi}.
It shows in particular that the niceness assumption
can be dropped for 
the implication \ref{i-Omtwo-exactly} $\imp$ \ref{i-Om}.

\begin{thm} \label{thm-reg-Phi-new}
Assume that $\bdy \Om \prec \bdytwo \Om$,
where $\bdy \Om$ is the given metric boundary.
Let $\Phi: \clOmtwo \to \clOm$ be the projection.
If
there are finitely many, and at least one, 
resolutive-irregular boundary points in $\Phi^{-1}(x_0)$
with respect to $\Omtwo$,
or more general there is a resolutive-irregular boundary point 
$\xh \in  \Phi^{-1}(x_0)$, with respect to $\Omtwo$,
which has arbitrarily small 
$\tautwo$-neighbourhoods $U$ such that 
\begin{equation} \label{eq-irr-cond}
   \bdytwo U \cap \{x' \in \Phi^{-1}(x_0): x' \text{ is 
resolutive-irregular with respect to $\Omtwo$}\}=\emptyset,
\end{equation}
then $x_0$ is irregular with respect to $\Om$.
\end{thm}

\begin{proof}
As $\xh$ is resolutive-irregular,
there is a resolutive $f \in C(\bdytwo \Om)$ such that
\[
     m:=\liminf_{\Om \ni y \totwo \xh} \Hpind{\Omtwo} f(y) <  f(\xh).
\]
We may assume that $f \ge 0$ and that $m < 1 < f(\xh)$.
Then there is a $\tautwo$-neighbourhood $U \subset \clOmtwo$ of
$\xh$ satisfying \eqref{eq-irr-cond} and
such that $f \ge 1$ on $\clUtwo \cap \bdytwo \Om$.
Let $G=U \cap \Om$ and
\[
       \ft=\begin{cases}
             1 & \text{on } \bdy G \cap \bdy \Om, \\
             \min\{1,\Hpind{\Omtwo} f\} & \text{on } \bdy G \cap\Om.
             \end{cases}
\]
Note that $\ft$ need not be continuous,
but since \eqref{eq-irr-cond} holds,
$\ft$ is continuous at $x_0$.

If $u \in \UU_{f}(\Omtwo)$, then $u \in \UU_{\ft}(G)$,
where $G$ is equipped with the given metric topology of $\clOm$.
Hence, 
$           \uHpind{G} \ft           \le \Hpind{\Omtwo} f$
in $G$.
It thus follows that
\[
   \liminf_{G \ni y \to x_0} \uHpind{G} \ft(y)
   \le    \liminf_{\Om \ni y \totwo \xh} \Hpind{\Omtwo} f(y)
    < 1=\ft(x_0),
\]
which, together with Proposition~\ref{prop-reg-1},
 shows that $x_0$ is irregular with respect to $G$.
Since $\bdy \Om$ is the given metric boundary,
Corollary~4.4 in Bj\"orn--Bj\"orn~\cite{BB}
(or Corollary~11.3 in \cite{BBbook})
shows that
$x_0$ is irregular with respect to $\Om$.
\end{proof}

\begin{proof}[Proof of Theorem~\ref{thm-reg-Phi}]
\ref{i-Om-reg} $\eqv$ \ref{i-Om} This
is a direct consequence of the fact that the given metric
boundary is resolutive.

\ref{i-Om}  $\imp$ \ref{i-Omtwo-at-least}
This follows from Proposition~\ref{prop-irreg-1=>2}.

\ref{i-Omtwo-at-least} $\imp$ \ref{i-Om-reg}
This follows from Theorem~\ref{thm-reg-Phi-new}.

\ref{i-Omtwo-at-least} $\imp$ \ref{i-Omtwo-exactly}
Let $\xh_1,\xh_2 \in \Phi^{-1}(x_0)$ be resolutive-irregular with respect to $\Omtwo$.
Assume that $\xh_1 \ne \xh_2$. 
We can proceed as in the proof of
Theorem~\ref{thm-reg-Phi-new}
finding functions $f_j$, $\ft_j$ and sets $U_j$, $G_j$ corresponding
to $\xh_j$, $j=1,2$. We may require that $U_1 \cap U_2 = \emptyset$.
As in the proof of Theorem~\ref{thm-reg-Phi-new}, we see that
$x_0$ is irregular with respect to $G_1$ and also with respect to $G_2$.
Since $G_1$ and $G_2$ are disjoint, this
contradicts 
Lemma~7.4 in Bj\"orn~\cite{ABcluster}
(or Lemma~11.32 in \cite{BBbook}).

\ref{i-Omtwo-exactly} $\imp$ \ref{i-Omtwo-at-least}
This is trivial.
\end{proof}

\begin{proof}[Proof of Theorem~\ref{thm-reg-bdym}]
Since $\bdym \Om$ is resolutive
(by   Theorem~8.2 in Bj\"orn--Bj\"orn--Shanmugalingam~\cite{BBSdir}),
and all points in $\bdy \Om$ split nicely with
respect to the Mazurkiewicz boundary $\bdym \Om$,
this follows directly from Theorem~\ref{thm-reg-Phi},
\end{proof}

A natural question is to which extent the assumptions
in  Theorem~\ref{thm-reg-Phi} are essential.
For two arbitrary compactifications $\bdyone \Om \prec \bdytwo \Om$,
the implication (corresponding to)
\ref{i-Om} $\imp$ \ref{i-Omtwo-at-least}
holds  by Proposition~\ref{prop-irreg-1=>2},
while \ref{i-Omtwo-exactly} $\imp$ \ref{i-Omtwo-at-least} is
trivial. 
Also obviously \ref{i-Om-reg} $\imp$ \ref{i-Om}, while the converse
implication fails by Example~\ref{ex-4-p>2}.
(All the counterexamples are collected in Section~\ref{sect-ex}.)
What about the other four implications not containing \ref{i-Om-reg}?

Consider first the case when $x_0$ splits nicely.
In this case Examples~\ref{ex-1} and~\ref{ex-3b} show
that no other implication holds.
In both cases $\clOmtwo=\clOm$ and both
 boundaries are resolutive (so that regularity and 
resolutive-regularity coincide).

If we instead keep the assumption that $\bdyone \Om=\bdy \Om$,
but drop the assumption that $x_0$ splits nicely,
then 
\ref{i-Omtwo-exactly} $\imp$ \ref{i-Om} 
by 
Theorem~\ref{thm-reg-Phi-new}.
On the other hand,
Example~\ref{ex-4} shows that 
\ref{i-Om} $\not\imp$ \ref{i-Omtwo-exactly}
and 
\ref{i-Omtwo-at-least} $\not\imp$ \ref{i-Omtwo-exactly},
even under the assumption that $\bdytwo \Om$ is resolutive.
We do not know if 
\ref{i-Omtwo-at-least}  $\imp$ \ref{i-Om}
holds without
the niceness assumption.

We also do not know if ``resolutive-regular'' can be replaced
by ``regular'' in Theorem~\ref{thm-reg-Phi}.
But if we also drop the niceness assumption, 
then 
\ref{i-Om} $\not\imp$ \ref{i-Omtwo-exactly},
\ref{i-Omtwo-at-least} $\not\imp$ \ref{i-Omtwo-exactly}
and 
\ref{i-Omtwo-at-least} $\not\imp$ \ref{i-Om},
by Examples~\ref{ex-4} and~\ref{ex-4-p>2},
while \ref{i-Om} $\imp$ \ref{i-Omtwo-at-least}
and \ref{i-Omtwo-exactly} $\imp$ \ref{i-Omtwo-at-least}
remain true,
and it is only the implication 
\ref{i-Omtwo-exactly} $\imp$ \ref{i-Om} that we do not know if it holds in
this case.

See Section~\ref{sect-further} for some sharper results
when $\Phi^{-1}(x_0)$ is finite.

\section{The Kellogg property and uniqueness results}
\label{sect-kellogg}

Our aim in this section is to establish 
Theorems~\ref{thm-Kellogg-bdym} and~\ref{thm-unique-bdym}
and suitable generalizations of them.
As an application of Theorem~\ref{thm-reg-Phi},
we can obtain the following so-called Kellogg property under the assumption
that $\bdy \Om \prec \bdytwo \Om$. 

\begin{thm} \label{thm-Kellogg}
\textup{(The resolutive Kellogg property)}
Assume that $\bdy \Om \prec \bdytwo \Om$,
where $\bdy \Om$ is the given metric boundary,
and that 
$\bCp(\,\cdot\,;\Om)$-q.e.\ 
boundary point in $\bdy \Om$ splits nicely.

Then $\bCp(\Irrres^2;\Omtwo)=0$,
where $\Irrres^2$
is the set of
resolutive-irregular boundary points with
respect to $\Omtwo$.
\end{thm}

If in addition 
$\bdytwo \Om$ is resolutive, then
we obtain the Kellogg property for $\bdytwo \Om$,
i.e.\ $\bCp(\Irr^2;\Omtwo)=0$, where 
$\Irr^2$ is the set of
irregular boundary points with
respect to $\Omtwo$.
Example~\ref{ex-3} shows that the Kellogg property
does not hold for arbitrary resolutive compactifications,
while Example~\ref{ex-4-p>2} shows that 
it does not hold for $\Irr^2$ with respect to 
arbitrary compactifications $\bdytwo \Om \succ \bdy \Om$.
Note also that by Proposition~\ref{prop-cap-all-bdy} the full boundary
of any compactification 
always has positive capacity, and hence the Kellogg property
is never trivial.

\begin{proof}
Let $\Irr \subset \bdy \Om$ be the set of irregular boundary points
with respect to $\Om$,
and $E \subset \bdy \Om$ be the set of boundary points which
do not split nicely.
Then $\Cp(\Irr)=0$, by the Kellogg property in Theorem~3.9 in 
Bj\"orn--Bj\"orn--Shan\-mu\-ga\-lin\-gam~\cite{BBS}
and Theorem~6.1 in Bj\"orn--Bj\"orn--Shan\-mu\-ga\-lin\-gam~\cite{BBS2}
(or Theorems~10.5 and~10.22 in \cite{BBbook}).
By Theorem~\ref{thm-reg-Phi}, $\Irrres^2 \subset \Phi^{-1}(\Irr \cup E)$,
where $\Phi: \clOmtwo \to \clOm$ is the projection.
Also, $\bCp(\Irr \cup E;\Om) \le \Cp(\Irr \cup E)=0$, by Lemma~5.2 in
Bj\"orn--Bj\"orn--Shanmugalingam~\cite{BBSdir}.
Hence, by Lemma~\ref{lem-bCp-Phi},
\[
     \bCp(\Irrres^2;\Om_1) 
     \le  \bCp(\Phi^{-1}(\Irr \cup E);\Om_1) = \bCp(\Irr \cup E;\Om) =0.
     \qedhere
\]
\end{proof}

\begin{proof}[Proof of Theorem~\ref{thm-Kellogg-bdym}]
Since $\bdym \Om$ is resolutive
(by   Theorem~8.2 in Bj\"orn--Bj\"orn--Shanmugalingam~\cite{BBSdir}),
and all points in $\bdy \Om$ split nicely with
respect to the Mazurkiewicz boundary $\bdym \Om$,
the Kellogg property for $\bdym \Om$
  follows directly from Theorem~\ref{thm-Kellogg}.
\end{proof}

To establish Theorem~\ref{thm-unique-bdym}, and its generalization 
Theorem~\ref{thm-unique} below,
we need the following two conditions:

\begin{enumerate}
\renewcommand{\theenumi}{\textup{(\roman{enumi})}}%
\item
  $\bdyone \Om$ is \emph{q.e.-invariant} if
  whenever $f \in C(\bdyone\Om)$ and 
  $h:\bdyone\Om\to\eR$ 
satisfies 
$
\bCp(\{x \in \bdyone \Om : h(x) \ne 0\};\Omone)=0,
$
then $\uHpind{\Omone} f = \uHpind{\Omone} (f+h)$;
\item
  the \emph{weak Kellogg property} holds for $\bdyone \Om$
  if for every $f \in C(\bdyone \Om)$ there is a set $E_f$, with
  $ \bCp(E_f;\Omone)=0$, such that
\begin{equation} \label{eq-weak-Kellogg}
              \lim_{\Om \ni y \toone x} \uHpind{\Omone} f(y) = f(x)
              \quad \text{for all } x \in \bdyone \Om \setm E_f.
\end{equation}
\end{enumerate}

One may think that we also need to require that $\bdyone \Om$ is resolutive,
but in fact this is a consequence of the two assumptions above, as
we show in Proposition~\ref{prop-inv+K=>res} below.
Note also that
by Proposition~\ref{prop-cap-all-bdy} the 
weak Kellogg property
is never trivial.
The equality in \eqref{eq-weak-Kellogg} can equivalently
be replaced by the inequality in Lemma~\ref{lem-res-reg},
see the proof of that lemma.
We do not know if all boundaries are q.e.-invariant.

\begin{prop} \label{prop-inv+K=>res}
  If $\bdyone \Om$ is q.e.-invariant and that the weak Kellogg
  property holds,
  then $\bdyone \Om$ is resolutive.
\end{prop}

Example~\ref{ex-4-p>2} shows that the weak Kellogg assumption cannot
be dropped, nor can it be replaced by the resolutive Kellogg property.

\begin{proof}
  Let $f \in C(\bdyone \Om)$
  and  $h=\infty \chi_{E_{-f}}$, where $E_{-f}$ comes
  from the weak Kellogg property for $-f$.
  Then $\lHpind{\Omone} f \in \UU_{f-h}$, and hence
  by the q.e.-invariance and Proposition~\ref{prop-lHp<=uHp},
  \[
     \uHpind{\Omone} f 
     = \uHpind{\Omone} (f-h)
     \le \lHpind{\Omone} f
     \le \uHpind{\Omone} f
     \]
from which it follows that $f$ is resolutive.     
\end{proof}

\begin{thm} \label{thm-unique}
  Assume that $\bdyone \Om$ is q.e.-invariant and that the weak Kellogg
  property holds.
Let $f \in C(\bdyone \Om)$. Then there exists
a unique bounded \p-harmonic function $u$ on $\Om$ such that
\begin{equation} \label{eq-unique}
    \lim_{\Om \ni y \toone x} u(y) = f(x) \quad 
\text{for $\bCp(\,\cdot\,;\Omone)$-q.e. } x \in \bdyone \Om.
\end{equation}
Moreover, $u=\Hpind{\Omone} f$.
\end{thm}

Examples~\ref{ex-3} and~\ref{ex-4-p>2} show that the weak Kellogg assumption cannot
be dropped, and the latter example also shows that it cannot
be replaced by the resolutive Kellogg property.

\begin{proof}
  By Proposition~\ref{prop-inv+K=>res}, $f$ is resolutive
  and, by the weak Kellogg property,
  $u=\Hpind{\Omone} f$ satisfies \eqref{eq-unique}, which
  establishes the existence.

  As for the uniqueness, let $u$ be a bounded \p-harmonic function
  and $E \subset \bdyone \Om$ be such that
  $\bCp(E;\Omone)=0$ and 
\[
    \lim_{\Om \ni y \toone x} u(y) = f(x) \quad 
\text{for } x \in \bdyone \Om \setm E.
\]
Let $h=\infty \chi_E$.
Then $u \in \UU_{f-h}$ and thus $u \ge \uHpind{\Omone} (f-h) = \Hpind{\Omone} f$,
by the q.e.-invariance.
By applying this to  $-u$ and $-f$ we also see that $u \le \Hpind{\Omone} f$.
Hence $u=\Hpind{\Omone} f$.
\end{proof}

\begin{proof}[Proof of Theorem~\ref{thm-unique-bdym}]
  By Theorem~8.2 in Bj\"orn--Bj\"orn--Shanmugalingam~\cite{BBSdir},
  $\bdym \Om$ is q.e.-invariant,
  and by Theorem~\ref{thm-Kellogg-bdym} the Kellogg property holds
  for $\bdym \Om$.
  Hence, the conclusion follows directly from Theorem~\ref{thm-unique}.
\end{proof}

\begin{cor} \label{cor-Sobolev-reg}
Assume that $\bdy \Om \prec \bdyone \Om$,
that $\bCp(\,\cdot\,;\Om)$-q.e.\ boundary point in $\bdy \Om$ splits nicely,
and that $\bdyone \Om$ is Sobolev-resolutive.
Then the assumptions, and therefore the conclusion,
in Theorem~\ref{thm-unique} hold     for $\bdyone \Om$.
\end{cor}

\begin{proof}
  As $\bdyone \Om$ is Sobolev-resolutive, it is automatically resolutive.
  Moreover it is q.e.-invariant by Proposition~7.1 in 
  Bj\"orn--Bj\"orn--Sj\"odin~\cite{BBSjodin}.
  It follows from Theorem~\ref{thm-Kellogg} that the Kellogg property
  holds.
  Hence the conclusion follows from Theorem~\ref{thm-unique}.
\end{proof}

The weak Kellogg property of course follows from the usual
Kellogg property (but not from the resolutive Kellogg property,
see Example~\ref{ex-4-p>2}).
However, if $\clOmone$ is metrizable then the weak Kellogg property
is equivalent to the usual Kellogg property.

\begin{thm}
Assume that $\clOmone$ is metrizable.
Then the weak Kellogg
  property holds if and only if the usual Kellogg property holds.
\end{thm}

Observe that we do not assume that $\clOmone$ is resolutive.
In the proof below we use that it follows from the metrizability
that $C(\bdyone \Om)$ is separable, and instead we could have used
this assumption. However, by
Theorem~2.10 in Bj\"orn--Bj\"orn--Sj\"odin~\cite{BBSjodin}
these two assumptions are equivalent.

The name ``weak Kellogg property'' was coined in
Bj\"orn~\cite{ABkellogg} where it was obtained
for quasiminimizers with respect to the given metric boundary.
For quasiminimizers, it is not known if the
Kellogg property holds or not.
Similarly, when
$C(\bdyone \Om)$ is not separable we do not know if the weak Kellogg
property for \p-harmonic functions implies the usual Kellogg property.

\begin{proof}
Assume that the weak Kellogg property holds.
  By Theorem~2.10 in \cite{BBSjodin},
  $C(\bdyone \Om)$ is separable, i.e.\ it contains
  a dense countable subset $A$.
  Let
  \[
     E = \bigcup_{f \in A} E_f,
  \]
     where $E_f$ comes from the weak Kellogg property for $f$.
     As the capacity is countably subadditive, we see that
     $\bCp(E;\Omone)=0$.

     If $f \in C(\bdyone \Om)$, then we can find $f_j \in A$ such
     that $f_j \to f$ uniformly.
     Then also $\uHpind{\Omone} f_j \to \uHpind{\Omone} f$ uniformly
and 
     it follows that
\[
   \lim_{\Om \ni y \toone x} \uHpind{\Omone} f(y) = f(x)
   \quad \text{for all } x \in \bdyone \Om \setm E.
\]
Hence $\Irr^1 \subset E$ and the Kellogg property follows.    
The converse implication is trivial.
\end{proof}

\section{Further results when \texorpdfstring{$\Phi^{-1}(x_0)$}{Phi-inverse(x0)}
is finite}
\label{sect-further}

The results in Section~\ref{sect-bdy-reg}--\ref{sect-kellogg} can be strengthened
when $\Phi^{-1}(x_0)$ is finite.
The following are the two main results in this section,
which improve upon  Proposition~\ref{prop-irreg-1=>2} and 
Theorem~\ref{thm-reg-Phi} under more restrictive assumptions.

\begin{prop} \label{prop-reg-singleton}
Assume that $\bdyone \Om \prec \bdytwo \Om$,
and let $\Phi: \clOmtwo \to \clOmone$ be the projection.
Let $x_0 \in \bdyone \Om$ 
and assume that $\Phi^{-1}(x_0)$ consists of just one point $\xh$.
Then $x_0$ is regular
with respect to $\Omone$
if and only if $\xh$ is regular 
with respect to $\Om^2$.
\end{prop}

Example~\ref{ex-3b} shows that Proposition~\ref{prop-reg-singleton}
cannot be generalized to the case when $\Phi^{-1}(x_0)$ is finite.
But if we assume that $\bdyone \Om = \bdy \Om$ is the given metric boundary,
we do obtain the following characterization.

\begin{thm} \label{thm-reg-Phi-finite}
Assume that $\bdy \Om \prec \bdytwo \Om$,
where $\bdy \Om$ is the given metric boundary.
Let $\Phi: \clOmtwo \to \clOm$ be the projection
and $x_0 \in \bdy \Om$.
Assume that $\Phi^{-1}(x_0)$ is finite.
Then the following are equivalent\/\textup{:}
\begin{enumerate}
\renewcommand{\theenumi}{\textup{(\alph{enumi}$'$)}}%
\item \label{j-Om}
The point $x_0$ is irregular with respect to $\Om$.
\item \label{j-Omtwo-at-least}
There is at least one irregular boundary point in $\Phi^{-1}(x_0)$
with respect to $\Omtwo$.
\item \label{j-Omtwo-exactly}
There is exactly  one irregular boundary point in $\Phi^{-1}(x_0)$
with respect to $\Omtwo$.
\end{enumerate}

Moreover, if $\xh \in \Phi^{-1}(x_0)$,
then $\xh$ is regular if and only if it is
resolutive-regular\/ \textup{(}with respect to $\clOmtwo$\textup{)}.
\end{thm}

The implication \ref{j-Omtwo-exactly} $\imp$ \ref{j-Omtwo-at-least} 
is of course trivial and the implication
\ref{j-Om} $\imp$ \ref{j-Omtwo-at-least}  follows from
Proposition~\ref{prop-irreg-1=>2}. One may ask how much of this result
remains true without assuming that the smaller boundary is the given metric
boundary
$\bdy \Om$. If we instead would assume that the larger boundary is  $\bdy \Om$,
then in fact no other than the two implications mentioned above would be
true, see Examples~\ref{ex-1} and~\ref{ex-3b}.

As a consequence of Theorem~\ref{thm-reg-Phi-finite}
we can obtain the non-resolutive Kellogg property under some conditions,
but without requiring resolutivity of the boundary.

\begin{thm} \label{thm-Kellogg-finite}
\textup{(The Kellogg property)}
Assume that $\bdy \Om \prec \bdytwo \Om$,
where $\bdy \Om$ is the given metric boundary,
and that 
$\Phi^{-1}{(x_0)}$ is finite for $\bCp(\,\cdot\,;\Om)$-q.e.\ $x_0 \in \bdy \Om$.

Then $\bCp(\Irr^2;\Omtwo)=0$.
\end{thm}

\begin{proof}
The proof is almost identical to the proof of 
Theorem~\ref{thm-Kellogg}, but using 
Theorem~\ref{thm-reg-Phi-finite}  instead of 
Theorem~\ref{thm-reg-Phi}.
\end{proof}

To prove Proposition~\ref{prop-reg-singleton}
we will use the following characterization,
which may be of independent interest.

\begin{prop} \label{prop-reg-char-fK}
Let $x_0 \in \bdyone \Om$.
  For each nonempty compact $K \subset \bdyone \Om \setm \{x_0\}$,
  let $f_K \in C(\bdyone \Om)$ be nonnegative and such that
  $f(x_0)=0 < \inf_K f_K$
  \textup{(}it exists by Tietze's extension theorem\/\textup{)}.

  Then $x_0$ is regular if and only if
  \begin{equation} \label{eq-reg-cond}
  \lim_{\Om \ni y \toone x_0} \uHpind{\Omone} f_K(y)=0
  \quad \text{for all nonempty compact }
  K \subset \bdyone \Om \setm \{x_0\}.
  \end{equation}
\end{prop}

If there is a nonnegative function $h \in C(\bdyone \Om)$
which is zero only at $x_0$, then we can use that function alone, i.e.\
we may let $f_K=h$ for all $K$.
This is however possible if and only if $x_0$ has a countable base of neighbourhoods
(which in particular holds if $\clOmone$ is first countable).

\begin{proof}
If $\bdyone \Om =\{x_0\}$ is the one-point compactification
of $\Om$, then $x_0$ is regular as it is the only boundary point,
and the equivalence is trivial. So assume that $\bdyone \Om \ne \{x_0\}$.

  Assume first that \eqref{eq-reg-cond} holds.
 Let  $f \in C(\bdyone \Om)$. We may assume that $f(x_0)=0$.
  Let $\eps >0$. Then there is a $\tauone$-neighbourhood $G \subsetneq \bdyone \Om$ 
of $x_0$
  such that $f < \eps $ in $G \cap \bdyone \Om$.
  Let $K=\bdyone \Om \setm G$ and $M=\sup_{\bdyone \Om} f / {\inf_K f_K}$.
  Then $f \le M f_K + \eps $ on $\bdyone \Om$.
  Hence
  \[
  \limsup_{\Om \ni y \toone x_0} \uHpind{\Omone} f(y)
  \le M    \lim_{\Om \ni y \toone x_0} \uHpind{\Omone} f_K(y) + \eps
  \le \eps.
  \]
  Letting $\eps \to 0$ shows that
  \[
  \limsup_{\Om \ni y \toone x_0} \uHpind{\Omone} f(y) \le 0 = f(x_0).
  \]
As $f$ was arbitrary,
Lemma~\ref{lem-res-reg} 
yields that $x_0$ is regular.
  The converse implication is trivial.  
\end{proof}

\begin{proof}[Proof of Proposition~\ref{prop-reg-singleton}]
One direction follows from Proposition~\ref{prop-reg-1}
but we will nevertheless show the full equivalence directly.

  For each nonempty compact $K \subset \bdyone \Om \setm \{x_0\}$,
  let $f_K \in C(\bdyone \Om)$ be nonnegative and such that
  $f(x_0)=0 < \inf_K f_K$
  (which exists by Tietze's extension theorem).
Also let $\fh_{\Kh} = f_{\Phi(\Kh)} \circ \Phi$
for each compact $\Kh \subset \bdytwo \Om \setm \{\xh\}$.

  If a net $\{y_\la\}_\la$ in $\Om$ converges to $x_0$ in $\clOmone$
  then it must converge to $\xh$ in $\clOmtwo$, and conversely.
  It thus follows from Theorem~\ref{thm-res-prec} that
  \[
  \lim_{\Om \ni y \toone x_0} \uHpind{\Omone} f_K(y)=0
  \quad \text{for all nonempty compact }
  K \subset \bdyone \Om \setm \{x_0\}
  \]
  if and only if
  \[
  \lim_{\Om \ni y \totwo \xh} \uHpind{\Omtwo} f_{\Kh}(y)=0
  \quad \text{for all nonempty compact }
  \Kh \subset \bdytwo \Om \setm \{\xh\},
  \]
  which together with Proposition~\ref{prop-reg-char-fK}
  completes the proof.
\end{proof}

When proving Theorem~\ref{thm-reg-Phi-finite}
we also need the following restriction result.

\begin{prop} \label{prop-reg-subset}
  Let $x_0 \in \bdyone \Om$,
  $G \subset \clOmone$ be a $\tauone$-neighbourhood of $x_0$
  and $\Uone=G \cap \Omone$.
  If $x_0$ is regular
  with respect to $\Uone$,
  then $x_0$ is regular with respect to $\Omone$.
\end{prop}

Boundary regularity with respect to the given metric is a local property
by Theorem~6.1 in Bj\"orn--Bj\"orn~\cite{BB} (or Theorem~11.11 in
  \cite{BBbook}).
Proposition~\ref{prop-reg-subset} shows that one direction
of this equivalence holds in full generality, while 
Example~\ref{ex-3b} shows that the other does not.
We will discuss this further in Section~\ref{sect-local-prop}.

\begin{proof}
  Let $f \in C(\bdyone \Om)$.
  By Tietze's extension theorem we can consider $f$ to be defined on $\clOmone$.
  We can also assume that $|f| \le 1$ and that $f(x_0)=0$.
  Let $K=\bdyone U \setm G$.
  Then, by Tietze's extension theorem, there is a nonnegative 
$h \in C(\bdyone U)$
  such that $h(x_0)=0$ and $h=1$ on $K$.
  Let now $\fh=\min\{h+f_\limplus,1\}$ on $\bdyone U$ and $\uh \in \UU_{\fh} (\Uone)$.
    Then
    \[
    u=\begin{cases}
    \min\{\uh,1\}, & \text{in } U, \\
    1, & \text{in } \Om \setm U
    \end{cases}
    \]
    is superharmonic in $\Om$, by Lemma~3.13
    in Bj\"orn--Bj\"orn--M\"ak\"al\"ainen--Parviainen~\cite{BBMP}
    (or \cite[Lemma~10.27]{BBbook}).
    It follows that $u \in \UU_f(\Omone)$
    from which we conclude that 
    $\uHpind{\Omone} f \le \uHpind{\Uone} \fh$ in $U$.
    Hence,
    \[
    \limsup_{\Om \ni y \toone x_0} \uHpind{\Omone} f(y)
    \le     \lim_{\Om \ni y \toone x_0} \uHpind{\Uone} \fh(y)
    = \fh(x_0) =f(x_0).
\]
As $f$ was arbitrary,
Lemma~\ref{lem-res-reg} 
yields that $x_0$ is regular with respect to $\Omone$.
\end{proof}

\begin{proof}[Proof of Theorem~\ref{thm-reg-Phi-finite}]
  Let $x_1,\ldots,x_m$ be the points in $\Phi^{-1}(x_0)$.
  As $\clOmtwo$ is normal we can for each $x_j$ find a
  $\tautwo$-neighbourhood $G_j$
  of $x_j$ whose closure avoids the other points.
  After having chosen all $G_j$ we can make each one smaller
  (if necessary, and still denoting it $G_j$)
  to make sure that its closure does not intersect
  the other closures either.

  Let next $U_j=G_j \cap \Om$.
  We will now consider $U_j$ both with its given metric
closure $\overline{U}_j$ and with the
  $\tautwo$-closure $\clUjtwo$.
  Note that $(\Phi|_\clUjtwo)^{-1}(x_0)=\{x_j\}$ and
  thus Proposition~\ref{prop-reg-singleton} is available.

  \ref{j-Om} $\imp$  \ref{j-Omtwo-at-least}
  This follows from Proposition~\ref{prop-irreg-1=>2}.

  $\neg$ \ref{j-Om} $\imp$ $\neg$ \ref{j-Omtwo-at-least}
  Let $j \in \{1,\ldots,m\}$.
  By Corollary~4.4 in Bj\"orn--Bj\"orn~\cite{BB} (or Corollary~11.3 in
  \cite{BBbook}) we see that $x_0$ is regular with respect to $U_j$.
  Hence, by Proposition~\ref{prop-reg-singleton},
  $x_j$ is regular with respect to $U_j^2$.
  Thus, $x_j$ is regular with respect to $\Om^2$,
  by Proposition~\ref{prop-reg-subset}.

  \ref{j-Omtwo-exactly} $\imp$ \ref{j-Omtwo-at-least}
  This is trivial.

  $\neg$ \ref{j-Omtwo-exactly} $\imp$ $\neg$ \ref{j-Omtwo-at-least}
  Assume first that there are (at least) two irregular boundary points
  in $\Phi^{-1}(x_0)$
  with respect to $\Omtwo$,
  which we may assume to be $x_1$ and $x_2$.
  It then follows from Proposition~\ref{prop-reg-subset}
  that $x_j$ is also irregular with respect to $U_j^2$, $j=1,2$.
  Hence, by Proposition~\ref{prop-reg-singleton},
  $x_0$ is irregular with respect to $U_j$, $j=1,2$.
Since $U_1$ and $U_2$ are disjoint, this
contradicts 
Lemma~7.4 in Bj\"orn~\cite{ABcluster}
(or Lemma~11.32 in \cite{BBbook}).
We thus conclude that if  \ref{j-Omtwo-exactly} fails, then
there is no irregular boundary point in 
$\Phi^{-1}(x_0)$ (with respect to $\Omtwo$),
and thus \ref{j-Omtwo-at-least} also fails.

The last part now follows from Theorem~\ref{thm-reg-Phi}
(together with the obvious fact that a regular point
is resolutive-regular).
\end{proof}

\section{Regularity as a local property}
\label{sect-local-prop}

Boundary regularity with respect to the given metric is a local property
by Theorem~6.1 in Bj\"orn--Bj\"orn~\cite{BB} (or Theorem~11.11 in
  \cite{BBbook}).
The following result extends this to a large class of boundaries
greater than the given metric boundary,
and also deduces a ``restriction result'' for the same boundaries
(quite different from the restriction result in Proposition~\ref{prop-reg-subset}).

\begin{thm} \label{thm-local-prop}
Assume that $\bdy \Om \prec \bdytwo \Om$,
where $\bdy \Om$ is the given metric boundary.
Let $\Phi: \clOmtwo \to \clOm$ be the projection,
$\xh \in \bdytwo \Om$ and $x_0 = \Phi(\xh)$.
Assume that either
\begin{enumerate}
\renewcommand{\theenumi}{\textup{(\roman{enumi})}}%
\item
$\Phi^{-1}(x_0)$ is finite\/\textup{;} or
\item \label{a-2}
$x_0$ splits nicely and $\bdytwo \Om'$ is resolutive for every open $\Om' \subset \Om$. 
\end{enumerate}

If $ U \subset \Om$ is an open set such that $\xh \in \clUtwo$
and 
$\xh$ is regular with respect to $\Omtwo$,
then $\xh$ is regular with respect to $\Utwo$.

Moreover, if 
$G^2 \subset \clOmtwo$ is a $\tautwo$-neighbourhood of $\xh$,
then $\xh$ is regular with respect to $\Omtwo$
if and only if it is regular with respect to  $(G \cap \Om)^2$.
\end{thm}

Example~\ref{ex-3b} shows that neither of these
two facts hold in general. 
We do not know if they may hold for arbitrary boundaries
larger than the given metric boundary.
As already noted, one direction in the equivalence does
hold for arbitrary compactifications, by Proposition~\ref{prop-reg-subset}.

\begin{proof}
We assume first that $\Phi^{-1}(x_0)$ is finite.
In order to prove the first part 
we need to consider two cases.

\medskip

\emph{Case}~1. $x_0$ is regular with respect to $\Om$ (with the
given metric boundary).
In this case it follows from 
Corollary~4.4 in Bj\"orn--Bj\"orn~\cite{BB} (or Corollary~11.3 in
  \cite{BBbook}), that $x_0$ is regular with respect to $U$.
It then follows from Theorem~\ref{thm-reg-Phi-finite} (applied
to $U$) that $\xh$ is regular with respect to $\Utwo$.

\medskip

\emph{Case}~1. $x_0$ is irregular with respect to $\Om$.
By Theorem~\ref{thm-reg-Phi-finite} there is
$x' \in \Phi^{-1}(x_0)$ which is irregular with respect to $\Omtwo$.
Since $\xh$ is regular  with respect to $\Omtwo$, we must have
$x' \ne \xh$.
As $\clOmtwo$ is a normal space there are $\tautwo$-neighbourhoods
$\Gh$ and $G'$ of $\xh$ and $x'$, respectively,
with disjoint $\tautwo$-closures.

Let $V=\Vh \cup V' :=(\Gh \cap U) \cup (G' \cap \Om)$.
By Proposition~\ref{prop-reg-subset}, $x'$ is irregular with respect to
$\Vtwo$.
Thus, by Theorem~\ref{thm-reg-Phi-finite}, 
\ref{j-Omtwo-at-least} $\imp$ \ref{j-Omtwo-exactly}, applied to $V$, 
$\xh$ must be regular with respect to $\Vtwo$.
As $\Vh$ and $V'$ are disjoint the Perron
solution $\uP_{\Vtwo} f$ within $\Vh$ only depends on the boundary values
on $\bdytwo \Vh$.
Since $\xh \notin \clVptwo$, it follows
that $\xh$ is regular also with respect to $\Vh^2$.
As $\Vh \subset U$, it follows from Proposition~\ref{prop-reg-subset}
that $\xh$ is regular
with respect to $\Utwo$.

\medskip

One direction of the second part follows directly from 
the first part, while the other one is a direct consequence 
of Proposition~\ref{prop-reg-subset}.

\medskip

The proof in case~\ref{a-2} is similar, but using
Theorem~\ref{thm-reg-Phi} instead of Theorem~\ref{thm-reg-Phi-finite}.
\end{proof}

The rather complicated condition \ref{a-2} above is essential for
our proof
(as Theorem~\ref{thm-reg-Phi} is applied to $V$)
 and it may seem hard to know when it is satisfied.
However, the main way of showing that the boundary $\bdytwo \Om$ is
resolutive (and almost the only available way) is to show that
continuous functions on $\bdytwo \Om$ can be uniformly approximated
(on $\bdytwo \Om$) by functions in
\[
    A=\{f: \clOmtwo \to \R : f \text{ is $\bCp(\,\cdot\,;\Omtwo)$-quasicontinuous 
       and } f\in \Np(\Om)\},
\]
from which the resolutivity (and even Sobolev-resolutivity)
of $\bdytwo \Om$ follows
by Theorem~6.4 and  Proposition~7.3 in Bj\"orn--Bj\"orn--Sj\"odin~\cite{BBSjodin}.
If one instead require that continuous functions on $\clOmtwo$ can be uniformly
approximated (on $\clOmtwo$) by functions in $A$, then not only
$\bdytwo \Om$ is resolutive, but also $\bdytwo \Om'$ for any open $\Om' \subset \Om$.
To see this one just need to take restrictions to $\clOmptwo$,
and apply the same resolutivity results.
Note however that it is not trivial that the restriction of a 
$\bCp(\,\cdot\,;\Omtwo)$-quasicontinuous function
is $\bCp(\,\cdot\,;\Omptwo)$-quasicontinuous,
since the conditions \eqref{eq-deff-bCp-1} and \eqref{eq-deff-bCp}
are different,
but this follows from the fact that $\bCp(\,\cdot\,;\Omtwo)$
is an outer capacity, by Proposition~4.2 in \cite{BBSjodin}.

In particular, this is true for the 
Mazurkiewicz boundary $\bdym \Om$,
if $\Om$ is 
finitely connected at the boundary,
since Lipschitz functions on $\clOmm$ belong to $\Np(\Om)$.
In this case, the resolutive of all subboundaries also follows
from Theorem~11.2 in Bj\"orn--Bj\"orn--Shanmugalingam~\cite{BBSdir}.

\begin{proof}[Proof of Theorem~\ref{thm-local-prop-bdym}]
It follows from 
Theorem~11.2 in Bj\"orn--Bj\"orn--Shan\-mu\-ga\-lin\-gam~\cite{BBSdir}
and the discussion after Definition~\ref{deff-nicely}
that condition \ref{a-2} in Theorem~\ref{thm-local-prop}
is satisfied, and thus the result follows
from Theorem~\ref{thm-local-prop}.
\end{proof}

\section{Counterexamples}
\label{sect-ex}

In this section we have collected a number of counterexamples
demonstrating the sharpness of our results (to the extent known to us).
These examples are all in $\R^2$.
To simplify notation we will consider $\R$ to be embedded into $\R^2$ in the
usual way.

The boundaries in Examples~\ref{ex-1}--\ref{ex-4} are all resolutive
(and thus obviously regularity and resolutive-regularity coincide),
while the boundary in Example~\ref{ex-4-p>2} is not.

\begin{example} \label{ex-1}
Let $\Om=B(0,2) \setm \{0,1\} \subset \R^2$ (unweighted) with $1 < p \le 2$.
Then $0$ and $1$ are both irregular with respect to $\R^2$.
Let $\clOmone$ be $\clOm$ with $0$ and $1$ identified,
$\hat{0}$ be the identified point, and
$\Phi : \clOm \to \clOmone$ be the projection.
It follows from Theorem~\ref{thm-res-prec} that $\bdyone \Om$ is resolutive.
Let $\fh=\chi_{\hat{0}}$ which is continuous on $\bdyone \Om$
and $f=\chi_{\{0,1\}}$ which is continuous on $\bdy \Om$.
By Theorem~\ref{thm-res-prec},
\[
     \Hpind{\Omone} \fh = \Hpind{\Om} f \equiv 0,
\]
which shows that $\hat{0}$ is irregular with respect to $\bdyone \Om$,
but $\Phi^{-1}(\hat{0})=\{0,1\}$ contains two irregular boundary
points.

This shows that even though $\hat{0}$ splits nicely  and $\bdyone \Om$
is resolutive,
\ref{i-Om} $\not\imp$ \ref{i-Omtwo-exactly}
and \ref{i-Omtwo-at-least}~$\not\imp$~\ref{i-Omtwo-exactly} 
in Theorem~\ref{thm-reg-Phi} 
for $\bdyone \Om \prec \bdy \Om$.
It also shows that 
\ref{j-Om} $\not\imp$ \ref{j-Omtwo-exactly}
and \ref{j-Omtwo-at-least}~$\not\imp$~\ref{j-Omtwo-exactly} 
in Theorem~\ref{thm-reg-Phi-finite} for $\bdyone \Om \prec \bdy \Om$.
\end{example}

\begin{example} \label{ex-3}
Let $p=3$,  $w(x)=|x|$ be a weight on $\R^2$, and $d\mu=w\,dx$.
Then $w$ is a Muckenhoupt $A_3$-weight, see Section~1.6
in Heinonen--Kilpel\"ainen--Martio~\cite{HeKiMa}.
Moreover 
$    \Cp(\{x\})>0$ if $x \ne 0$, while $\Cp(\{0\})=0$,
see Example~2.22 in \cite{HeKiMa}.

Let $\Om$, $\clOmone$, $f$, $\fh$ and $\Phi$ be as in Example~\ref{ex-1}.
It again follows from Theorem~\ref{thm-res-prec} that $\bdyone \Om$ is resolutive.
This time
\[
    \bCp(\{\hat{0}\};\Omone)=\bCp(\{0,1\};\Om) 
    \ge \bCp(\{1\};\Om) >0,
\]
by Lemma~\ref{lem-bCp-Phi} and a straightforward
calculation.
Moreover
\[
   \Hpind{\Omone} \fh = \Hpind{\Om} f = \Hpind{\Om \cup \{0\}} f=:u,
\]
by Theorem~\ref{thm-res-prec} and 
 Bj\"orn--Bj\"orn--Shanmugalingam~\cite[Corollary~6.2]{BBS2}
(or  \cite[Corollary~10.22]{BBbook}).

The set $\Om \cup \{0\}$ is regular by the Kellogg property
(for the given metric boundary).
Hence $\lim_{x \to 1} u(x)=1$ and $\lim_{x \to 2} u(x)=0$,
showing that $u$ is nonconstant. 
Thus  by the strong maximum principle, $0<u<1$ in $\Om \cup \{0\}$.
In particular 
\[
    \liminf_{\Om \ni y \toone \hat{0}} u(y) = u(0)<1,
\]
showing that $\hat{0}$ is irregular.
Hence both the weak and the usual Kellogg properties fail for $\Omone$,
and also the conclusion in Theorem~\ref{thm-unique} fails even
though $\bdyone \Om$ is q.e.-invariant (as the empty set is the only
boundary
set with zero capacity).
\end{example}

\begin{example} \label{ex-3b}
Let $p=3$, $X=[0,8] \times [0,1] \subset \R^2$ and equip it
with the measure $d\mu= w\, dx$, where $w(x)=|x|$.
Also let $\Om=X \setm \{0,8\}$.
As in Example~\ref{ex-3} we have 
$\Cp(\{0\})=0$ and  
$\Cp(X \setm \Om)=\Cp(\{8\})>0$.
Then $0$ is irregular, 
while $8$ is regular with respect to $\Om$.
Let $\clOmone$ be $\clOm$ with $0$ and $8$ identified as 
$\hat{0}$, and
$\Phi : \clOm \to \clOmone$ be the projection.
It follows from Theorem~\ref{thm-res-prec} that $\bdyone \Om$ is resolutive.
Moreover,  $\hat{0}$ is regular with respect to $\Omone$ (as it is
the sole boundary point), while $\Phi^{-1}(\{\hat{0}\})$ contains exactly
one
irregular boundary point. 

This shows that even though $\hat{0}$ splits nicely  and $\bdyone \Om$
is resolutive,
\ref{i-Omtwo-exactly} $\not \imp$ \ref{i-Om}
and \ref{i-Omtwo-at-least} $\not \imp$ \ref{i-Om} 
in Theorem~\ref{thm-reg-Phi}
and 
\ref{j-Omtwo-exactly} $\not \imp$ \ref{j-Om}
and \ref{j-Omtwo-at-least} $\not \imp$ \ref{j-Om} 
in Theorem~\ref{thm-reg-Phi-finite} for $\bdyone \Om \prec \bdy \Om$.
It also shows that Proposition~\ref{prop-reg-singleton}
cannot be generalized to the case when $\Phi^{-1}(x_0)$ is finite.

Let 
$
G = (\Om \setm B(4,2)) \cup \{\hat{0}\}
$,
which is a $\tauone$-neighbourhood of $\hat{0}$,
$U=\{x \in G \cap \Om : |x| < 4\}$ and 
$f=\chi_{\bdyone U \setm \{\hat{0}\}} \in C( \bdyone G)$.
Since $0$ is irregular with respect to $G$, we see
that $\Hpind{G^1} f = \chi_U$
and thus $\hat{0}$ is irregular with respect to $(G \cap \Om)^1$.
Hence the converse implication to the one in Proposition~\ref{prop-reg-subset}
does not hold in general,
and regularity is not a local property in this situation.
In particular, neither of the two parts in Theorem~\ref{thm-local-prop}
hold in this case.
\end{example}

\begin{example} \label{ex-4}
Let $\Om = B(0,1) \setm \{0\} \subset \R^2$ (unweighted).
We want to create a compactification of $\Om$ such that the
outer boundary $\bdy \Om \setm \{0\}$ is as for the given metric,
but such that the boundary point $0$ corresponds to
an interval $I=[-1,1]$.
To do so, we let $h(x)=\sin (1/|x|)$.
Let also $h_j(x_1,x_2)=x_j$ be the coordinate functions,
and $Q=\{h,h_1,h_2\}$.

We want $\clOmtwo$ to be the smallest compactification such
that the three functions $h$, $h_1$ and $h_2$ have continuous
extensions to $\clOmtwo$.
Such an extension exists 
and  is called 
the
\emph{$Q$-compactification} $\clOmQ$ of $\Om$, see 
Bj\"orn--Bj\"orn--Sj\"odin~\cite{BBSjodin},
(it is unique up to homeomorphism).
Note that $\clOm \prec \clOmtwo$, and that
$\clOmtwo$ is metrizable, by Theorem~2.10 in \cite{BBSjodin}.
Let $\Phi : \clOmtwo \to \clOm$ be the corresponding projection.
A neighbourhood base for $a \in I$ is given by 
\begin{equation} \label{eq-nbhd-base}
    U_\eps :=\{x \in I : |x-a| < \eps\} 
    \cup \{x \in \Om : |x| < \eps \text{ and }
    |h(x)-a|<\eps\},
\end{equation}
with $0 < \eps <1$.

Assume that $1 < p \le 2$.
Then $0$ is irregular with respect to $\Om$, 
and $\uHpind{\Om} f\equiv 0$ for  $f=\chi_{\{0\}}$.
Then, by Theorem~\ref{thm-res-prec},
$P_{\Omtwo}({f \circ \Phi})= P_{\Om}f \equiv 0$, showing
that all points in  $I=\Phi^{-1}(\{0\})$ are irregular.
Moreover, in the terminology of Bj\"orn~\cite{ABclass}
(or \cite[Chapter~13]{BBbook}),
$0$ is semiregular and the Perron solutions with respect
to $\Om$ ignore the value at $0$ for continuous boundary data.
Hence, if $f \in C(\bdytwo \Om)$ and we let
\begin{equation} \label{eq-f1-f2}
   f_1(x)=\begin{cases}
      f(x), & x \in \bdy \Om \setm \{0\}, \\
      \inf_{ I} f, & x=0,
      \end{cases}
   \quad \text{and} \quad 
   f_2(x)=\begin{cases}
      f(x), & x \in \bdy \Om \setm \{0\}, \\
      \sup_{I} f, & x=0,
      \end{cases}
\end{equation}
then, by Theorem~\ref{thm-res-prec},
\[
  \uHpind{\Omtwo} f 
   \le \Hpind{\Omtwo} (f_2 \circ \Phi)
  =  \Hpind{\Om} f_2 
  =  \Hpind{\Om} f_1 
   =  \Hpind{\Omtwo} (f_1 \circ \Phi)
   \le \lHpind{\Omtwo} f,
\]
from which we conclude that $f$ is resolutive and thus
$\bdytwo \Om$ is resolutive.
This shows that 
\ref{i-Om} $\not \imp$ \ref{i-Omtwo-exactly}
and \ref{i-Omtwo-at-least} $\not \imp$ \ref{i-Omtwo-exactly}
in Theorem~\ref{thm-reg-Phi} 
when $x_0$ does not split nicely, 
even if $\bdytwo \Om$ is assumed to be resolutive.
\end{example}

\begin{example} \label{ex-4-p>2}
Assume now that $\Om$ and $\clOmtwo$ are as in
Example~\ref{ex-4}, but
this time with  $p>2$. 
In this case $\Cp(\{0\})>0$ and thus $0$ is regular with 
respect to $\Om$.
Let again $f \in C(\bdytwo \Om)$ and let $f_1$ and $f_2$ be given
by \eqref{eq-f1-f2}.
If $a < \sup_I f$ then, because of \eqref{eq-nbhd-base},
 any function $u \in \UU_f(\Omtwo)$ 
is necessarily greater than $a$ on concentric 
circles which are arbitrarily close to $0$, and hence
it is greater than $a$ in a neighbourhood
of $0$, 
by the minimum principle for superharmonic 
functions, see  Heinonen--Kilpel\"ainen--Martio~\cite[Theorem~7.12]{HeKiMa} 
(or  \cite[Theorem~9.13]{BBbook}).
As this holds for all $a < \sup_I f$,
we see that $u \in \UU_{f_2}(\Om)$.
Conversely, any function in $\UU_{f_2}(\Om)$ necessarily 
belongs to $\UU_f(\Omtwo)$.
We therefore conclude that
\[
   \uHpind{\Omtwo} f = \Hpind{\Om} f_2.
\] 
Similarly, $\lHpind{\Omtwo} f = \Hpind{\Om} f_1$.
As $\Hpind{\Om} f_1 \equiv \Hpind{\Om} f_2$ if and only if $f$ is constant
on $I$, only such $f$ are resolutive with respect to $\Omtwo$,
and thus $\bdytwo \Om$ is not resolutive.

From this we can easily conclude that all the points
in $I$ are irregular, but resolutive-regular.
As 
$\bCp(I;\Omtwo) = \bCp(\{0\},\Om)>0$, by Lemma~\ref{lem-bCp-Phi},
we see that neither the weak nor the usual Kellogg property hold 
with respect to $\Omtwo$.
On the other hand the resolutive Kellogg property does hold,
as there are no resolutive-irregular
boundary points.
Moreover, this shows that 
\ref{i-Omtwo-at-least} $\imp$ \ref{i-Om}
in Theorem~\ref{thm-reg-Phi} can fail
when $x_0$ does not split nicely, 
even if $\bdytwo \Om$ is assumed to be resolutive.

In fact, a similar argument (using concentric circles) shows that
\begin{equation} \label{eq-E}
\bCp(E;\Omtwo) = \bCp(\{0\},\Om)>0
\quad \text{for all nonempty sets } E \subset I.
\end{equation}
Hence $\bdytwo \Om$ is q.e.-invariant (as the only set
$E \subset \bdytwo \Om$ with zero capacity is the empty set).
Since any bounded \p-harmonic function  on $\Om$ has a limit
as $x \to 0$ (see below),  we also  conclude from \eqref{eq-E} that if
$f \in C(\bdytwo \Om)$ is nonconstant on $I$, then 
the (existence) conclusion in Theorem~\ref{thm-unique} fails for $f$.
We also see that the weak Kellogg property in
Proposition~\ref{prop-inv+K=>res} and Theorem~\ref{thm-unique}
neither can be dropped nor replaced by the resolutive Kellogg property.

It remains to show that
for any  bounded  \p-harmonic function $u \ge 0$ on $\Om$
the limit $\lim_{x \to 0} u(x)$ exists.
To this end, let 
\[
m(r)=\inf_{|x|=r} u(x)
 \text{ and } 
 M(r)=\sup_{|x|=r} u(x),
 \quad 0 <r<1,
\]
which are both continuous functions that, by
the strong maximum principle 
(see  \cite[Theorem~7.12]{HeKiMa} or  \cite[Theorem~9.13]{BBbook}),
can have at most one local extreme point each.
Hence the limits $m:=\lim_{r \to 0} m(r)$ and $M:=\lim_{r \to 0} M(r)$ exist.
By Harnack's inequality, there is a constant $A$ such that
$M\bigl(\tfrac12\bigr) \le Am\bigl(\tfrac12\bigr)$,
and by scaling invariance we can apply it also in smaller punctured balls.
Let $\eps >0$. Then there is $\rho >0$ such that for $0<r<\rho$ we have
$u > m -\eps$ in $B(0,2r)$.
Applying the Harnack inequality to $u-(m-\eps)$ shows that
\[
    M(r)-m+\eps \le A(m(r) - m+\eps).
\]
Letting first $r \to 0$ and then $\eps \to 0$ shows that $M=m$.     

\end{example}

\end{document}


\authortitle{Anders Bj\"orn}
{Correction of ``The Kellogg property  for \p-harmonic functions ...''}
\title{Correction of ``The Kellogg property and boundary regularity 
for \p-harmonic functions 
with respect to the Mazurkiewicz boundary and other compactifications''}
\author{
Anders Bj\"orn \\
\it\small Department of Mathematics, Link\"oping University, \\
\it\small SE-581 83 Link\"oping, Sweden\/{\rm ;}
\it \small anders.bjorn@liu.se
}

\date{}
\maketitle

\noindent{\small
{\bf Abstract.} 
We fill in a gap in the proofs of Theorems~1.1--1.4 in
``The Kellogg property and boundary regularity for \p-harmonic functions
with respect to the Mazurkiewicz boundary and other compactifications'',
to appear in \emph{Complex Var.\ Elliptic Equ.},
doi:10.1080/17476933.2017.1410799.
} 

\bigskip
\noindent
    {\small \emph{Key words and phrases}:
boundary regularity,
Mazurkiewicz boundary,
split nicely.
}

\medskip
\noindent
{\small Mathematics Subject Classification (2010): 
Primary: 31C45; Secondary: 31E05, 35J66, 35J92, 49Q20.
}

\bigskip

It has come to my attention that there is a gap in the argument
showing that 
every boundary point $x_0 \in  \bdy \Om$
splits nicely with respect to the Mazurkiewicz boundary
if
$\Om$ is as in Theorems~1.1--1.4 in \cite{ABgenreg}
(i.e.\ $\Om \subset X$ is a bounded domain which is finitely connected at the boundary,
where $X$ is a complete metric space equipped with a doubling
measure supporting a \p-Poincar\'e inequality, $1<p<\infty$).
This fact is mentioned after Definition~6.1, tacitly assuming that 
$\Phi^{-1}(x_0)$ is at most countable.
However, $\Phi^{-1}(x_0)$ can be uncountable 
even under the assumptions in Theorems~1.1--1.4 in \cite{ABgenreg},
see Example~7.5 in Bj\"orn--Bj\"orn--Li~\cite{BBLi}.
Nevertheless, $x_0$ does split nicely in this case.

To see this, let $\xh \in \Phi^{-1}(x_0)$ and let 
$V$ be a Mazurkiewicz neighbourhood of $\xh$.
(The Mazurkiewicz metric is always defined with respect to $\Om$.)
Then there is $r>0$ so 
that $V \supset B^M(\xh,3r):=\{x \in \clOmm: d_M(x,\xh)<3r \}$.
Let $G$ be the component of $\Om \cap B\bigl(x,r\bigr)$ which has $\xh$
in its Mazurkiewicz closure,
and let 
\[
U=\clGm \setm \Phi^{-1}(\{x \in \clOm : d(x,x_0)=r\}),
\]
which is an open  Mazurkiewicz neighbourhood of $\xh$.
As $G$ is connected,
we see that 
$d_M(x,y) \le \diam G \le 2r$ whenever $x,y \in G$.
Thus
\[
     B^M(\xh,r) \subset U \subset B^M(\xh,3r) \subset V.
\]
Moreover, if $x' \in \clUm \setm U$, then $d(\Phi(x'),x_0)=r$
and in particular $\Phi(x') \ne x_0$.
Hence $x_0$ splits nicely.

\begin{ack}
The author was supported by the Swedish Research Council,
grant 2016-03424.
\end{ack}